\documentclass[onefignum,onetabnum]{siamart190516}

\usepackage{todonotes}

\usepackage{amsfonts}
\usepackage{graphicx}
\usepackage{epstopdf}
\usepackage{algorithmic}
\usepackage{nicefrac}
\usepackage{stmaryrd}
\ifpdf
  \DeclareGraphicsExtensions{.eps,.pdf,.png,.jpg}
\else
  \DeclareGraphicsExtensions{.eps}
\fi
\usepackage{booktabs}

\newsiamremark{remark}{Remark}
\newsiamremark{hypothesis}{Hypothesis}
\crefname{hypothesis}{Hypothesis}{Hypotheses}
\newsiamthm{claim}{Claim}

\headers{Multilevel Spectral Domain Decomposition}{P. Bastian, R. Scheichl, L. Seelinger and A. Strehlow}

\title{Multilevel Spectral Domain Decomposition
}

\author{Peter Bastian\thanks{Interdisciplinary Center for Scientific Computing (IWR), Universität Heidelberg, 
  (\email{peter.bastian@iwr.uni-heidelberg.de}, \url{https://conan.iwr.uni-heidelberg.de/people/peter/}).}
\and Robert Scheichl\thanks{Department of Applied Mathematics and IWR, Universität Heidelberg,
  (\email{r.scheichl@uni-heidelberg.de},\email{linus.seelinger@iwr.uni-heidelberg.de},\email{arne.strehlow@uni-heidelberg.de}).}
\and Linus Seelinger\footnotemark[3]
\and Arne Strehlow\footnotemark[3]}

\usepackage{amsopn}

\DeclareMathOperator{\supp}{supp}
\DeclareMathOperator{\mathspan}{span}
\DeclareMathOperator{\range}{range}

\ifpdf
\hypersetup{
  pdftitle={Multilevel Spectral Domain Decomposition},
  pdfauthor={P. Bastian, R. Scheichl, L. Seelinger and A. Strehlow}
}
\fi

\newcommand{\new}[1]{{{#1}}}

\usepackage[normalem]{ulem}

\begin{document}

\maketitle

\begin{abstract}
Highly heterogeneous, anisotropic coefficients, e.g. in the simulation
of carbon-fibre composite components, can lead to extremely
challenging finite element systems. Direct solvers for the resulting
large and sparse linear systems suffer from severe memory requirements
and limited parallel scalability, while iterative solvers in general lack robustness. 
Two-level spectral domain decomposition methods can provide such
robustness for symmetric positive definite linear systems, by using coarse spaces
based on independent generalized eigenproblems in the subdomains. Rigorous
condition number bounds are independent of mesh size, number of subdomains, as
well as coefficient contrast. However, their parallel scalability is still
limited by the fact that (in order to guarantee robustness) the coarse
problem is solved via a direct method. In this paper, we introduce a multilevel variant in
the context of subspace correction methods and provide a general convergence
theory for its robust convergence for abstract, elliptic variational
problems. Assumptions of the theory are verified for
conforming, as well as for discontinuous Galerkin methods applied to a scalar diffusion problem.
Numerical results illustrate the performance of the method for two- and three-dimensional
problems and for various discretization schemes, in the context of scalar diffusion and linear elasticity.
\end{abstract}



\section{Introduction}

In this paper we are concerned with the solution of 
large and sparse linear systems
\begin{equation}\label{eq:linear_system}
Ax = b
\end{equation}
where $A\in\mathbb{R}^{n \times n}$ is a symmetric and positive definite (SPD) matrix
arising from the discretization of an elliptic (system of) partial
differential equation(s) (PDE) on a bounded domain $\Omega
  \subset \mathbb{R}^d$, with $d=2,3$ the spatial dimension.

Direct methods for solving \eqref{eq:linear_system} are very effective for relatively small problems
but suffer from severe memory requirements
and limited parallel scalability \cite{lexcellent:tel-00737751,REINARZ2018269}.
We focus on iterative methods instead. In PDE applications, the mesh parameter $h$
needs to be chosen sufficiently small to control the error in the numerical solution, leading to very
large systems $n\sim h^{-d}$. A variety of methods have been
developed in the past that converge (almost) independently of the mesh
size $h$ as well as of the number of
processors $p \sim H^{-d}$ -- assuming a parallel partitioning of
$\Omega$ into subdomains with diameter bounded by $H$. These include
in particular multigrid (MG) and domain decomposition (DD) methods 
\cite{Smith96,TosWid05,DoleanDD}. These methods require a ``coarse solver component''  
providing global information transfer. In this paper we focus on extensions of the two-level overlapping
Schwarz DD method.

While robustness with respect to (w.r.t.) $h$ and $H$ is achieved by many methods,
robustness w.r.t.~mesh anisotropy or to large variations
in problem parameters, such as the permeability coefficient in porous media flow or the Lam\'e parameters
in linear elasticity, are more difficult to achieve. Construction of robust coarse
spaces was a theme in multigrid early on \cite{Alcouffe1981} and lead to the development of
algebraic multigrid (AMG) methods \cite{xu_zikatanov_2017}. Rigorous convergence bounds for
aggregation-based AMG applied to nonsingular symmetric M-matrices with nonnegative row sums
are provided in \cite{Notay2012}. 
In the context of overlapping DD methods, it was shown in
\cite{SVZ2012}  based on weighted Poincaré inequalities 
\cite{Pechstein2013} that the standard method can be robust
w.r.t.~strong coefficient variation inside subdomains, but this puts hard constraints on the domain
decomposition. Construction of coarse spaces based on multiscale basis functions \cite{Aarnes2005,GLS2007}
can be very effective, but also leads to no rigorous robustness w.r.t.~arbitrary
coefficient variations.
A breakthrough was achieved by constructing coarse spaces based on solving certain local generalized eigenvalue
problems (GEVP). In \cite{NATAF20101163}, a local eigenvalue problem
involving the Dirichlet-to-Neumann map
was introduced, and later analysed in \cite{DtN2012} based on \cite{Pechstein2013}.
Different GEVP were proposed in \cite{Effendiev2010,Effendiev2012,Geneo14} together
with an analysis that is solely based on approximation properties of the chosen local eigenfunctions.

In this paper, we consider extensions of the GenEO ({\bf Gen}eralized
{\bf E}igenproblems in the {\bf O}verlap) method
introduced in \cite{Geneo14}. The method works on general elliptic systems of PDEs and is rigorously shown to
be robust w.r.t.~coefficient variations when all eigenfunctions w.r.t.~eigenvalues below a threshold are
included into the coarse space. This number can be related to the number of isolated high-conductivity regions,
\cite{Effendiev2010}, but also depends on the specific GEVP used. The local GEVP can either be constructed
from local stiffness matrices or in a fully algebraic way from the global stiffness matrix
through a symmetric positive semi-definite (SPSD) splitting \cite{GenEOAlgebraic}.
An extension of the GenEO approach to nonoverlapping domain decomposition methods
as well as non-selfadjoint problems can be found in \cite{HAFERSSAS2015959}.

Two-level domain decomposition methods traditionally employ direct solvers in the subdomains
and on the coarse level. While iterative solvers could be employed, they then need to be robust with
respect to coefficient variations as well. If such a solver would be available, it could be used instead
of the domain decomposition method. Thus we assume that such solver is not available.
The use of direct solvers in the subdomains and on the coarse level puts an upper limit on the size
of the subdomain/coarse problem and thus on the total problem size due to run-time and memory requirements.
The more severe penalty is typically imposed by the memory
requirements, in particular for todays supercomputers which have many cores with relatively little
memory per core. In order to give a concrete example, consider the HAWK system of HLRS, Stuttgart, Germany.
It consists of 5632 nodes, each containing 2 CPUs with 64 cores each and 256 GB of memory, i.e. 2GB of memory per core.
This limits the fine level subdomain size to about $n_0=250000$ degrees of freedom (dof) in scalar, three-dimensional problems.
In order to maintain a good speedup (in the setup phase), the coarse system size is then also limited to 250000 dof
or about 12500 subdomains (or cores) when we assume 20 basis vectors per subdomain from the GEVP.
Thus the GenEO method would not scale to the full Hawk machine.
This estimate could be improved by employing parallel direct solvers, but the scalability of these
solvers is limited \cite{lexcellent:tel-00737751,REINARZ2018269,GenEO19}.

Scalability beyond $10^4$ cores can be achieved by employing more than two levels.
A robust multilevel method employing spectral coarse spaces based on the work \cite{Effendiev2010}
has been proposed in \cite{Willems2014}. It employs a hierarchy of finite element meshes, GEVPs that
are very similar to the ones that we propose in this paper and a nonlinear AMLI-cycle. 
Robustness and level independence is achieved with a $W$-cycle, which is not desirable
in a parallel method. While robustness with respect to coefficient variations is proven and 
demonstrated numerically, the problem sizes are rather small.
A multilevel version of GenEO was proposed in \cite{GenEO19} based on the
SPSD splitting introduced in \cite{GenEOAlgebraic}. 

In this paper, we formulate a natural extension of GenEO \cite{Geneo14} to an
arbitrary number of levels. In contrast to \cite{GenEO19}, our method is formulated in a variational framework based on
subspace correction \cite{Xu92}. The convergence theory of \cite{Geneo14} is
generalized to nonconforming discretizations as well as to multiple levels and several different GEVPs.
This enables us to prove robust convergence  also for discontinuous
Galerkin methods suitable for high coefficient contrast \cite{Ern01042009}.
Condition number bounds are derived for an iterated two-level method as in \cite{GenEO19} but
also for a fully additive multilevel preconditioner. 
Numerical results up to $2^{16}$ subdomains and more than $10^8$ dofs demonstrate the effectiveness
of the approach.

The paper is organized as follows. In \cref{sec:main} we formulate
the multilevel spectral domain decomposition method in a variational framework.
In \cref{sec:theory} we provide the analysis of the method and briefly describe
our implementation in \cref{sec:impl}. Numerical results follow
in Section  \ref{sec:results} and we end with conclusions in
\cref{sec:conclusions}.

\section{Multilevel Spectral Domain Decomposition}
\label{sec:main}

\subsection{Subspace Correction} 

Throughout the paper we assume the linear system \eqref{eq:linear_system} arises
by inserting a basis representation into the variational problem
\begin{equation}\label{eq:variational_problem}
u_h \in V_h: \qquad a_h(u,v) = l_h(v) \qquad \forall v \in V_h,
\end{equation}
where $V_h$ is a finite element space (a finite dimensional vector space), $a_h : V_h \times V_h \to \mathbb{R}$
is a symmetric positive definite bilinear form and $l_h\in V_h^\prime$ is a linear form.
The subscript $0<h\in\mathbb{R}$ denotes that $a_h$ and $l_h$ are defined on a shape regular mesh 
$\mathcal{T}_h$ consisting of elements $\tau$ with diameter at most $h$, partitioning the
domain $\Omega\subset\mathbb{R}^d$. Elements $\tau=\mu_\tau(\hat\tau)$ are assumed to
be open and the image of a reference element $\hat\tau$ under the diffeomorphism $\mu_\tau$.
Reference elements are either the reference simplex or the reference cube in dimension $d$.

Subspace correction methods \cite{Xu92} are based on a splitting
\begin{equation*}
V_h = V_{h,1} + \ldots + V_{h,p}
\end{equation*}
of $V_h$ into $p$ possibly overlapping subspaces $V_{h,i}$. Any such splitting gives rise to the
iterative parallel subspace correction method
\begin{gather}
u_h^{k+1} = u_h^{k} + \omega \sum_{i=1}^p w_i^k, \quad\text{with $w_i^k$ given by} \label{eq:subspace_iteration}\\
w_i^k\in V_{h,i}: \qquad  a_h(w_i^k,v) = l(v)-a_h(u_h^{k},v) \qquad \forall v\in V_{h,i} . \nonumber
\end{gather}
Here, $\omega>0$ is a suitably chosen damping factor. Sequential subspace correction, see \cite{Xu92},
typically converges faster but offers less parallelism. Parallel subspace correction is also called
additive subspace correction and sequential subspace correction is called multiplicative subspace correction due
to the form of the error propagation operator. Hybrid variants may offer a good compromise in practice.

Practical implementation of subspace correction employs a
basis representation
\begin{align*}
V_{h} &= \text{span} \{\phi_1, \ldots, \phi_n \}, & V_{h,i} &= \text{span} \{\phi_{i,1}, \ldots, \phi_{i,n_i} \}, &
\phi_{i,j} &= \sum_{k=1}^{n} \left(R_i \right)_{j,k} \phi_k.
\end{align*}
The rectangular matrices $R_i$ represent the basis of $V_{h,i}$ in terms of the basis of $V_h$.
Expanding $u_h^k = \sum_{j=1}^{n} \left( x^k \right)_j \phi_j$ leads to the algebraic form
\begin{align*}
x^{k+1} &= x^k + \omega \sum_{i=1}^p R_i^T A_i^{-1} R_i (b-Ax^{k})
\end{align*}
with $(A)_{r,s} = a_h(\phi_s, \phi_r)$,  $A_i = R_i A R_i^T$ and the preconditioner
$ B = \sum_{i=1}^p R_i^T A_i^{-1} R_i $. $B$ is typically used as a preconditioner in
the conjugate gradient method.

Multilevel spectral domain decomposition methods are introduced below in the framework
of subspace correction. They employ a decomposition
\begin{equation}\label{eq:mlsdd}
V_h = \sum_{l=0}^L V_{h,l} = V_{h,0} + \sum_{l=1}^L \sum_{i=1}^{P_l} V_{h,l,i}
\end{equation}
where $L+1$ is the number of levels with $L$ being the finest level and $0$ the coarsest
level. Identifying $V_h = V_{h,L}$ we have the nested level-wise spaces $V_{h,l} \subset V_{h,l+1}$
for $0 \leq l < L$. On each level $l>0$, $V_{h,l}$ is split into $P_{l}$ subspaces $V_{h,l,i}$.

\subsection{Hierarchical Domain Decomposition}

\begin{figure}[tbp]
\begin{center}
\includegraphics[width=0.33\textwidth]{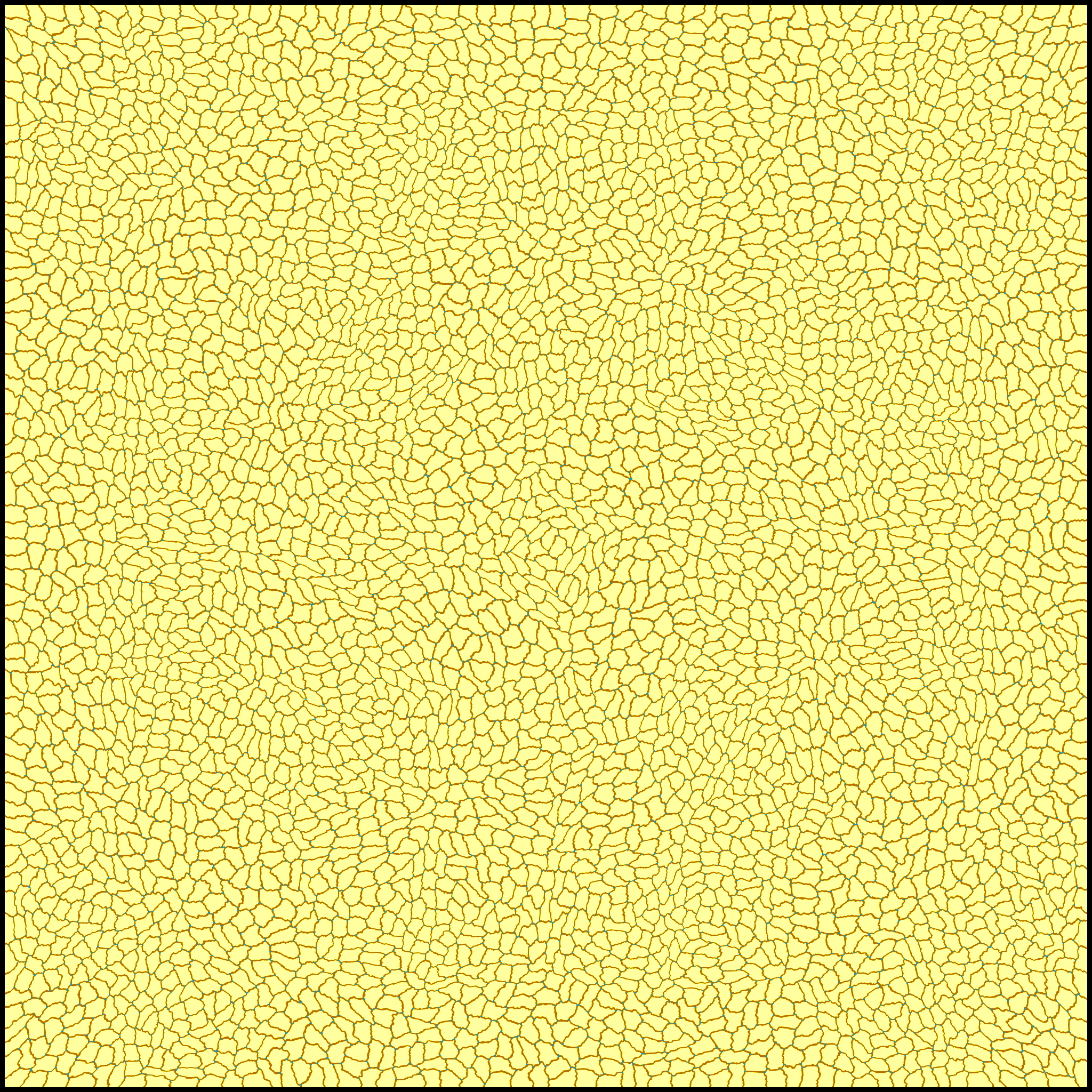}\hfill
\includegraphics[width=0.33\textwidth]{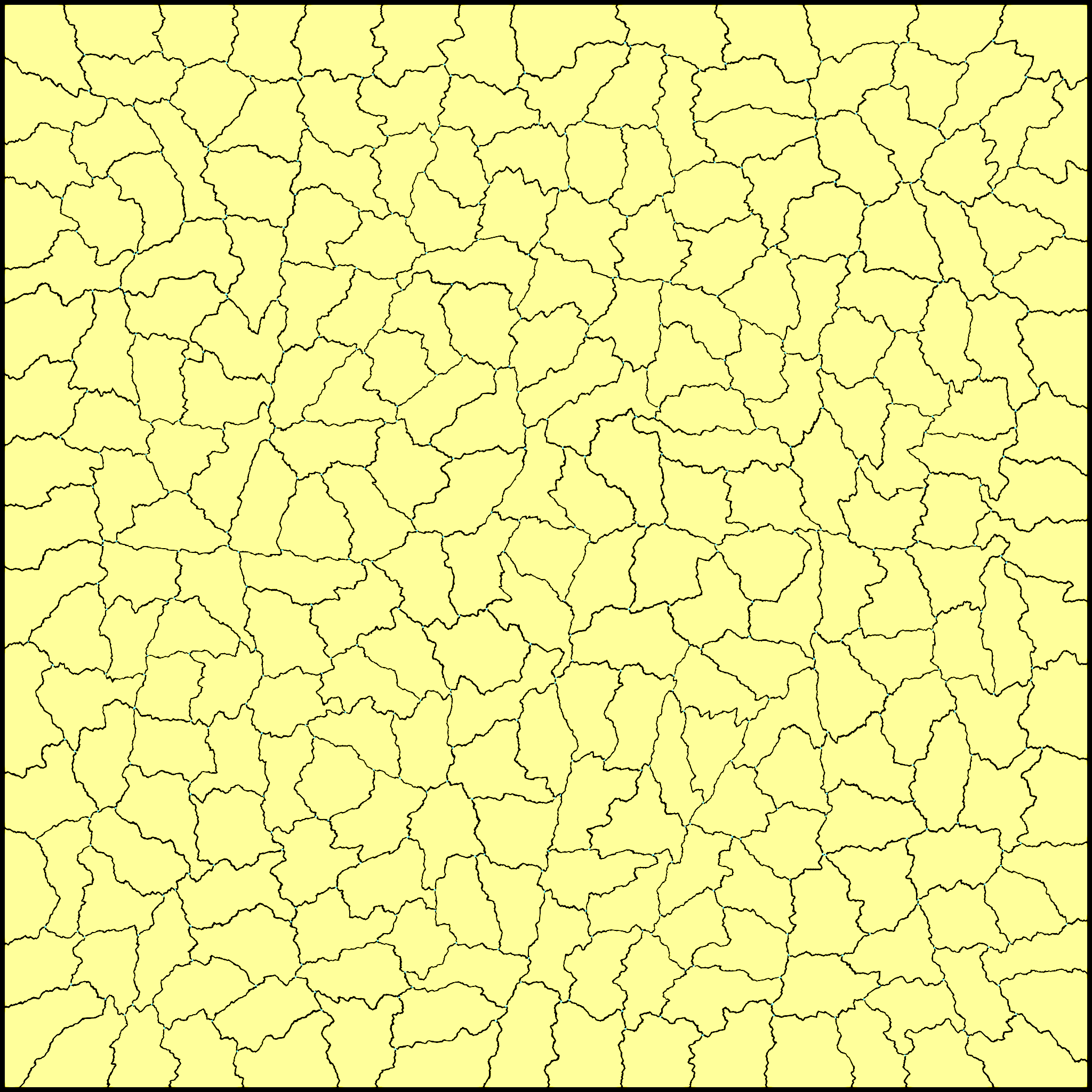}\hfill
\includegraphics[width=0.33\textwidth]{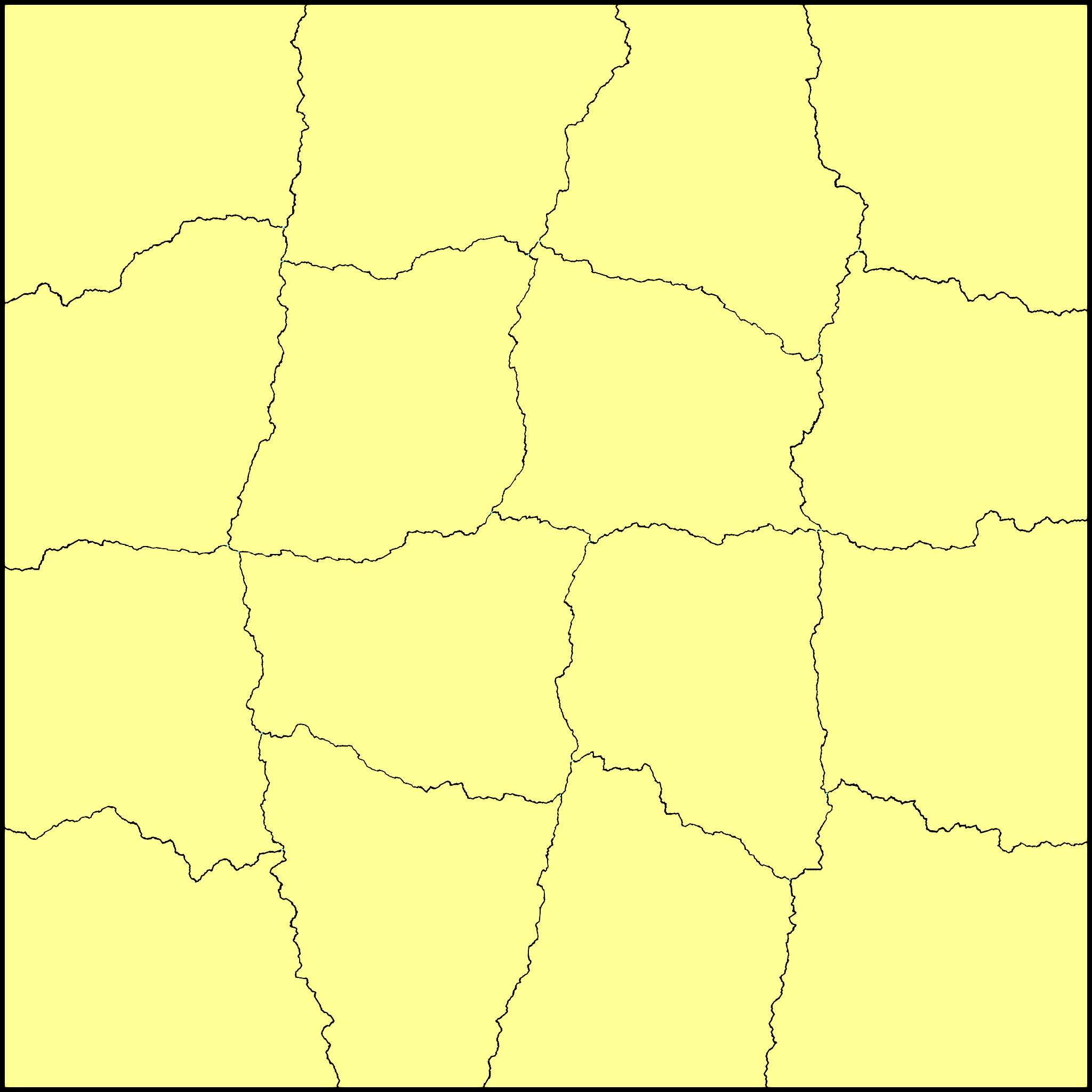}
\end{center}
\caption{Multilevel domain decomposition in two dimensions. A triangular mesh with 1.6 million vertices
is partitioned into 4096 (level 3), 256 (level 2) and 16 (level 1) nested subdomains using ParMetis \cite{KARYPIS199896}. 
No overlap is used on the coarse levels.}
\label{fig:mldd_2d}
\end{figure}

The construction of the subspaces is related to a hierarchical decomposition of the domain $\Omega$
into subdomains $\Omega_{l,i}$ for $0< l \leq L$ and $ 1\leq i \leq P_{l}$.
``Hierarchical'' means that each subdomain $\Omega_{l,i}$ is the union of subdomains
on the next finer level $l+1$. In particular, our multilevel method
employs a single fine mesh $\mathcal{T}_h$ given by the user.
All subdomains are unions of elements of the mesh $\mathcal{T}_h$. Figure \ref{fig:mldd_2d} shows a 
decomposition of a triangular mesh
employing three levels.

The domain decomposition is obtained from a decomposition of
$\mathcal{T}_h$ as follows:
\begin{enumerate}
\item On the finest level $L$, decompose $\mathcal{T}_h$ into $P_{L}$ overlapping sets 
\begin{align*}
\mathcal{T}_{h,L,i} &\subset \mathcal{T}_h, \quad 1 \leq i \leq P_{L}, 
& \bigcup_{i=1}^{P_{L}} \mathcal{T}_{h,L,i}  &= \mathcal{T}_h ,
\end{align*}
by first partitioning $\mathcal{T}_h$ into $P_{L}$ nonoverlapping sets of elements using
a graph partitioner such es ParMetis \cite{KARYPIS199896} and then adding a user defined
overlap $\delta$ in terms of layers of elements.
\item On levels $0< l < L$, determine decompositions of the subdomain index sets
\begin{align*}
J_{l,i}&\subset \{1,\ldots,P_{l+1}\}, \quad 1 \leq i \leq P_{l},
& \bigcup_{i=1}^{P_{l}} J_{l,i} &= \{1,\ldots,P_{l+1}\} .
\end{align*}
The sets $J_{l,i}$ may overlap but are not required to. Such decomposition is
obtained by a graph partitioner using the subdomain graph instead of the mesh.
With this we obtain the mesh partitioning on level $l$ as
\begin{equation*}
\mathcal{T}_{h,l,i} = \bigcup_{k\in J_{l,i}} \mathcal{T}_{h,l+1,k} .
\end{equation*}
\item The subdomains $\Omega_{l,i}$ are now defined from the mesh decomposition as
\begin{align*}
\Omega_{l,i} &= \text{Interior}\biggl( \bigcup_{\tau\in\mathcal{T}_{h,l,i}} \bar\tau \biggr),
&& 0<l \leq L, \quad 1 \leq i \leq P_{l}. 
\end{align*}
\item On the coarsest level 0 we employ always only one subdomain.
\end{enumerate}
\pagebreak


For finite volume and discontinuous Galerkin methods we need to introduce notation for
mesh faces. $\gamma$ is an interior face if it is the
intersection of two elements $\tau^-(\gamma), \tau^+(\gamma)\in\mathcal{T}_h$ and  has dimension $d-1$.
All interior faces are collected in the set $\mathcal{F}_h^I$. Likewise, 
$\gamma$ is a boundary face if it is
the intersection of some element $\tau^-(\gamma)\in\mathcal{T}_h$ with $\partial\Omega$ of dimension $d-1$.
All boundary faces make up the set $\mathcal{F}_h^{\partial\Omega}$.
With each $\gamma\in\mathcal{F}_h^I$ we associate
a unit normal vector $\nu_\gamma$ oriented from $\tau^{-}(\gamma)$ to $\tau^+(\gamma)$.
For a boundary face $\gamma\in\mathcal{F}_h^{\partial\Omega}$ 
its unit normal $\nu_\gamma$ coincides with the unit normal to $\partial\Omega$.
Related to the submeshes $\mathcal{T}_{h,l,i}$ we define the corresponding sets of faces
\begin{subequations}\label{eq:DD_faces}
\begin{align}
\mathcal{F}_{h,l,i}^{I} &= \{ \gamma\in\mathcal{F}_h^I : \tau^-(\gamma)\in\mathcal{T}_{h,l,i} \wedge \tau^+(\gamma)\in\mathcal{T}_{h,l,i}\}, \\
\mathcal{F}_{h,l,i}^{\partial\Omega} &= \{ \gamma\in\mathcal{F}_h^{\partial\Omega} : \tau^-(\gamma)\in\mathcal{T}_{h,l,i}\} .
\end{align}
\end{subequations}

\subsection{Spectral Coarse Space Construction}

Based on the hierarchical domain decomposition we can now formulate the construction of the
coarse spaces $V_{h,l}\subseteq V_h$ and $V_{h,l,i}$ introduced in \eqref{eq:mlsdd}.
First define the auxiliary local spaces
\begin{equation} \label{eq:def_neumann_space}
\overline V_{h,l,i} = \{ v|_{\Omega_{l,i}} : v \in V_{h,l} \} .
\end{equation}
The restriction operator $r_{l,i} : V_{h,l} \to \overline V_{h,l,i}$, $r_{l,i} v = v|_{\Omega_{l,i}}$, restricts
the domain of a finite element function. For $v$ being zero on $\partial\Omega_{l,i}\cap\Omega$,
the extension operator $e_{l,i} : \overline V_{h,l,i} \to V_{h,l}$ extends functions by zero outside of $\Omega_{l,i}$.
Another major ingredient is the partition of unity.

\begin{definition}
A partition of unity on level $0< l \leq L$ is a family of operators 
$\chi_{l,i} : \overline V_{h,l,i} \to \overline V_{h,l,i}$ such that
\begin{enumerate}
\item $\chi_{l,i} v$ is zero on $\partial\Omega_{l,i}\cap\Omega$  for all $v \in \overline V_{h,l,i}$ and
\item $\sum_{i=1}^{P_{l}} e_{l,i}\chi_{l,i}r_{l,i} v = v$ for all $v \in V_{h,l}$,
\end{enumerate}
\end{definition}

Restriction and extension operators as well as the partition of unity operators with the required properties can be defined
for conforming finite element spaces as well as discontinuous Galerkin finite element spaces.

\noindent The construction of the coarse spaces is now recursive over the levels:
\begin{enumerate}
\item Set $V_{h,L} = V_h$. Set $l=L$.
\item For each subdomain $1 \leq i \leq P_{l}$ solve a GEVP
\begin{equation} \label{eq:theGEVP}
w_{l,i,k} \in \overline{V}_{l,i} : \qquad \overline a_{l,i}(w_{l,i,k},v) = \lambda_{l,i,k} \overline b_{l,i}(w_{l,i,k},v) \qquad \forall v \in \overline V_{l,i}
\end{equation} 
with appropriately defined local bilinear forms $\overline a_{l,i}$ and $\overline b_{l,i}$ detailed below.
Let eigenvalues and corresponding eigenvectors be ordered s.t. $\lambda_{l,i,k} \leq \lambda_{l,i,k+1}$.
\item With a user defined parameter $\eta$ define the coarse space $V_{h,l-1}$ as
\begin{equation} \label{eq:coarse_space_def}
V_{h,l-1} = \bigoplus_{i=1}^{P_l} \mathspan\left\{ \phi_{l,i,k} = e_{l,i} \chi_{l,i} w_{l,i,k} : \lambda_{l,i,k}<\eta \right\} .
\end{equation}
\item Set $l=l-1$. If $l>0$ goto step 2, otherwise stop.
\end{enumerate}
The subspaces for the subspace correction method based on $V_{h,l}$ are then given by
\begin{equation}\label{eq:dirichlet_space}
V_{h,l,i} =\left\{ v  \in V_{h,l} : \supp v \subset \overline{\Omega}_{l,i} \right\} \subset \overline V_{h,l,i}, \qquad 0<l\leq L, \quad 1\leq i \leq P_l
\end{equation}
together with $V_{h,0}$.

\begin{remark} The most important properties of this construction are:
\begin{itemize}
\item Within each level all GEVPs can be solved in parallel.
\item Basis functions $\phi_{l,i,k}\in V_h$ have support only in $\Omega_{l,i}$.
\item Functions in $V_{h,l,i}$ are zero at their subdomain boundary (except on the global Neumann boundary).
Functions in $\overline V_{h,l,i}$ are not necessarily zero at their subdomain boundary (except on the global Dirichlet boundary).
\end{itemize}
\end{remark}

\section{Convergence Theory}
\label{sec:theory}

\subsection{Standard Subspace Correction Theory}

For the standard analysis of subspace correction methods, the following
$a_h$-orthogonal projections $\mathcal{P}_{l} : V_h \to V_{h,l}$ and 
$\mathcal{P}_{l,i} : V_h \to V_{h,l,i}$ are introduced \cite{Xu92,TosWid05}:
\begin{align*}
a_h(\mathcal{P}_{l} w,v) &= a_h(w,v),  \quad \forall v \in V_{h,l}, &
a_h(\mathcal{P}_{l,i} w,v) &= a_h(w,v),  \quad \forall v \in V_{h,l,i}.
\end{align*}
This allows to write the error propagation 
operator $\mathcal{E}$ of the iteration \eqref{eq:subspace_iteration} as
\begin{align*}
\mathcal{E} &= \mathcal{I} - \omega \mathcal{P}, & \mathcal{P} = \mathcal{P}_0 + \sum_{l=1}^{L} \sum_{i=1}^{P_l} \mathcal{P}_{l,i}.
\end{align*}
Taking $\omega=1/\lambda_{\max}(\mathcal{P})$,
the convergence factor of the iteration \eqref{eq:subspace_iteration} is $\rho = 1-1/\kappa_2(\mathcal{P})$
with the spectral condition number
$\kappa_2(\mathcal{P}) = \lambda_{\max}(\mathcal{P})/\lambda_{\min}(\mathcal{P})$.
The goal of the analysis below is to provide upper and lower bounds of the form
\begin{equation*}
\gamma a_h(v,v) \leq a_h(\mathcal{P}v,v) \leq \Gamma a_h(v,v) \qquad \forall v \in V_h,
\end{equation*}
which in turn give a bound on the condition number $\kappa_2(\mathcal{P}) \leq \Gamma/\gamma$.

The analysis is based on two major properties.

\begin{definition}[Coloring and domain decomposition] \label{def:coloring}

\begin{enumerate}
\item We say the multilevel domain decomposition admits a level-wise coloring with $k_0\in\mathbb{N}$ colors, if
for each level $0<l\leq L$ there exists a map $c_l : \{1,\ldots,P_l\} \to \{1,\ldots,k_0\}$ such
that $$i\neq j \wedge c_l(i)=c_l(j) \Rightarrow a_h(v_i,v_j) = 0 \qquad \forall v_i \in V_{h,l,i}, v_j \in V_{h,l,j}.$$
\item The multilevel domain decomposition is called admissible, if on each level $0<l\leq L$ every interior face $\gamma$ 
is interior to at least one subdomain.
\end{enumerate}
\end{definition}

\begin{lemma} \label{lem:upper_bound}
Let the multilevel domain decomposition have a finite coloring with $k_0$ colors.
Then parallel subspace correction satisfies the upper bound
\begin{equation*}
a_h(\mathcal{P}v,v) \leq (1+k_0 L)  a_h(v,v) \qquad \forall v \in V_h.
\end{equation*}
\vspace{-6mm}
\begin{proof}
See \cite[p. 182]{Smith96}.
\end{proof}
\end{lemma}

\begin{definition}[Stable splitting] \label{def:stable_splitting}
The subspaces $V_{h,0}$ and $V_{h,l,i}$, $0<l\leq L$, $1 \leq i \leq P_l$, admit a stable splitting
if there exists $C_0>0$ and for each $v\in V_h$ a
decomposition $v = v_0 + \sum_{i=1}^L \sum_{i=1}^{P_l} v_{l,i}$, $v_0 \in V_{h,0}$, $v_{l,i} \in V_{h,l,i}$, such that
\begin{equation*}
a_h(v_0,v_0) + \sum_{i=1}^L \sum_{i=1}^{P_l} a_h(v_{l,i},v_{l,i}) \leq C_0 \, a_h(v,v) .
\end{equation*}
\end{definition}

\begin{lemma} \label{lem:lower_bound}
If the subspaces $V_{h,0}$ and $V_{h,l,i}$, $0<l\leq L$, $1 \leq i \leq P_l$ admit
a stable splitting with constant $C_0$, parallel subspace correction satisfies the lower bound
\begin{equation*}
C_0^{-1}  a_h(v,v)  \leq a_h(\mathcal{P}v,v) \qquad \forall v \in V_h.
\end{equation*}
\vspace{-6mm}
\begin{proof}
See e.g. \cite{TosWid05}.
\end{proof}
\end{lemma}

\begin{theorem} \label{thm:k2bound}
If the subspaces $V_{h,0}$ and $V_{h,l,i}$, $0<l\leq L$, $1 \leq i \leq P_l$ have a finite coloring with $k_0$ colors and admit
a stable splitting with constant $C_0$, parallel subspace correction satisfies the bound
\begin{equation*}
\kappa_2(\mathcal{P})  \leq C_0 (1+k_0 L) .
\end{equation*}
\vspace{-6mm}
\begin{proof}
Follows immediately from Lemma \ref{lem:upper_bound} and Lemma \ref{lem:lower_bound}.
\end{proof}
\end{theorem}

\subsection{Abstract Schwarz Theory for Spectral Domain Decomposition}

We now generalize the theory in \cite{Geneo14} to the multilevel case
and to more general discretization schemes.
For each application, the following three definitions have to be verified.

\begin{definition}[Strengthened triangle inequality under the square] \label{def:binomial}
The domain decomposition allows a strengthened triangle inequality under the square if
there exists $a_0>0$ independent of $0<l\leq L$ and $P_l$ such that
for any collection of $v_{l,i} \in V_{h,l,i}$:
\begin{equation*}
\left\| \sum_{i=1}^{P_l} v_{l,i} \right\|_{a_h}^2 \leq a_0 \sum_{i=1}^{P_l} \left\| v_{l,i} \right\|_{a_h}^2.
\end{equation*}
\end{definition}
Here $\|v\|_{a_h} = \sqrt{a_h(v,v)}$ is the norm induced by $a_h$.

\begin{definition}[Positive semi-definite splitting] \label{def:posdefsplitting}
The bilinear forms $\overline a_{l,i}$ introduced in \eqref{eq:theGEVP} provide a positive semi-definite
splitting of $a_h$ if there exists $b_0>0$ independent of $0<l\leq L$ and $P_l$ such that for each level $0<l\leq L$:
\begin{equation*}
\sum_{i=1}^{P_l} \left| r_{l,i} v_l \right|_{\overline a_{l,i}}^2 \leq b_0 \left\|v_l \right\|_{a_h}^2, \qquad \forall v_l \in V_{h,l} .
\end{equation*} 
\end{definition}
Here $|v|_{\overline a_{l,i}} = \sqrt{\overline a_{l,i}(v,v)}$ is the semi-norm induced by $\overline a_{l,i}$.
Symmetric positive semi-definite splittings on the algebraic level were introduced in \cite{GenEOAlgebraic,GenEO19}. 

\begin{definition}[Local stability] \label{def:local_stability}
The subspaces $V_{h,l}$ and $V_{h,l,i}$ are called locally stable if
there exists a constant $C_1>0$ idependent of $l$ and for each $v_l \in V_{h,l}$, $1 \leq l \leq L$, a
decomposition $v_l = v_{l-1} + \sum_{i=1}^{P_l} v_{l,i}$ with $v_{l-1}\in V_{h,l-1}$ and
$v_{l,i}\in V_{h,l,i}$  such that the following inequalities hold:
\begin{equation*}
\left\| v_{l,i} \right\|_{a_h}^2 \leq C_1 \left| r_{l,i} v_{l} \right|_{\overline a_{l,i}}^2, \qquad 1 \leq i \leq P_l .
\end{equation*}
\end{definition}

\begin{lemma}[Levelwise stability] \label{lem:alternative_stability}
Let the spectral coarse spaces admit strengthened triangle inequalities under the square, 
let the bilinear forms $\overline a_{l,i}$ provide a symmetric positive semi-definite splitting
and let the subspaces $V_{h,l,i}$ be locally stable.
Then, for each level $0<l\leq L$ there exists a decomposition $v_l = v_{l-1} + \sum_{i=1}^{P_l} v_{l,i}$ 
with $v_{l-1}\in V_{h,l-1}$ and $v_{l,i}\in V_{h,l,i}$ such that
\begin{align*}
\left\| v_{l-1}\right \|_{a_h}^2 \leq 2(1+ a_0 b_0 C_1) \left\|  v_{l}  \right \|_{a_h}^2
\end{align*}
and 
\begin{align*}
\left\| v_{l-1}\right \|_{a_h}^2 + \sum_{i=1}^{P_l} \left\| v_{l,i} \right\|_{a_h}^2  
\leq (2 +b_0C_1(1+2a_0)) \left\| v_{l} \right\|_{a_h}^2 .
\end{align*}
\vspace{-3mm}
\begin{proof} 
Adaptation of \cite[Lemma 2.9]{Geneo14}.
From Def. \ref{def:local_stability} and Def. \ref{def:posdefsplitting} we get
\begin{align*}
\sum_{i=1}^{P_l} \left\| v_{l,i} \right\|_{a_h}^2 \leq C_1\sum_{i=1}^{Pl}  \left| r_{l,i} v_l \right|_{\overline a_{l,i}}^2 
\leq b_0 C_1 \left\| v_l \right\|_{a_h}^2 .
\end{align*}
Thus, it follows from $v_{l-1} = v_l -  \sum_{i=1}^{P_l}
v_{l,i}$ and from Def. \ref{def:binomial} that
\begin{align*}
\left\| v_{l-1}\right \|_{a_h}^2 &= \left\|  v_l -  \sum_{i=1}^{P_l} v_{l,i}  \right \|_{a_h}^2
\leq 2 \left\|  v_l  \right \|_{a_h}^2 + 2 \left\|  \sum_{i=1}^{P_l} v_{l,i}  \right \|_{a_h}^2
\leq 2 \left\|  v_l  \right \|_{a_h}^2 + 2 a_0 \sum_{i=1}^{P_l} \left\| v_{l,i} \right\|_{a_h}^2\\
&\leq 2 \left\|  v_l  \right \|_{a_h}^2  + 2 a_0 b_0 C_1 \left\| v_l \right\|_{a_h}^2
= 2(1+ a_0 b_0 C_1) \left\|  v_l  \right \|_{a_h}^2 .
\end{align*}
Combining both results yields
\begin{align*}
\left\| v_{l-1}\right \|_{a_h}^2 + \sum_{i=1}^{P_l} \left\| v_{l,i} \right\|_{a_h}^2  
&\leq (2 +b_0C_1(1+2a_0)) \left\| v_l \right\|_{a_h}^2 .
\end{align*}
\end{proof}
\end{lemma}

The two-level decomposition from Definition \ref{def:local_stability} implies a multilevel decomposition
as follows: For any given function $v_h\in V_h$:
\begin{equation} \label{eq:multilevel_decomposition}
V_h \ni v_h = v_L = v_0 + \sum_{l=1}^L \sum_{i=1}^{P_l} v_{l,i}, \qquad \text{with} \qquad v_l -v_{l-1} = \sum_{i=1}^{P_l} v_{l,i} .
\end{equation}
We can now formulate the first major result of our paper.
\begin{lemma}[Multilevel stability] \label{Lem:multilevel_stability}
Let the spectral coarse spaces admit strengthened triangle inequalities under the square, 
let the bilinear forms $\overline a_{l,i}$ provide a symmetric positive semi-definite splitting
and let the subspaces $V_{h,l,i}$ be locally stable.
Then the multilevel decomposition \eqref{eq:multilevel_decomposition} satisfies for any $v_h\in V_h$
\begin{equation*}
\|v_0\|_{a_h}^2  + \sum_{l=1}^L \sum_{i=1}^{P_l} \left\| v_{l,i} \right\|_{a_h}^2 \leq 
C^{L} \left( 1 + \frac{b_0 C_1 }{C-1} \right) \| v_h \|_{a_h}^2
\end{equation*}
with $C=2(1+ a_0 b_0 C_1)$.
\begin{proof} For any level $0<l\leq L$, we obtain as in the proof of Lemma \ref{lem:alternative_stability}
\begin{align*}
\sum_{i=1}^{P_l} \left\| v_{l,i} \right\|_{a_h}^2 \leq C_1\sum_{i=1}^{P_l}  \left| r_{l,i} v_l \right|_{\overline a_{l,i}}^2 
\leq b_0 C_1 \left\| v_l \right\|_{a_h}^2 .
\end{align*}
Furthermore, Lemma \ref{lem:alternative_stability} implies
\begin{align*}
\|v_l\|_{a_h}^2 \leq \left[ 2(1+ a_0 b_0 C_1)\right]^{L-l} \| v_h \|_{a_h}^2
\end{align*}
for $0\leq l \leq L$. Together we obtain with setting $C=2(1+ a_0 b_0 C_1)$
\begin{align*}
\|v_0\|_{a_h}^2  &+ \sum_{l=1}^L \sum_{i=1}^{P_l}  \left\| v_{l,i} \right\|_{a_h}^2 
\leq C^{L} \| v_h \|_{a_h}^2 +  b_0 C_1  \| v_h \|_{a_h}^2 \sum_{l=1}^L  C^{L-l}\\
&\leq C^{L} \| v_h \|_{a_h}^2 +  b_0 C_1 \frac{C^L}{C-1} \| v_h \|_{a_h}^2
= C^{L} \left( 1 + \frac{b_0 C_1 }{C-1} \right) \| v_h \|_{a_h}^2.
\end{align*}
\end{proof}
\end{lemma}

\begin{theorem}\label{thm:multilevelrate}
Let the assumptions of Lemma \ref{Lem:multilevel_stability} hold and let the domain decomposition
admit a coloring with $k_0$ colors on each level.
Then the parallel multilevel subspace correction method satisfies the condition number bound
\begin{equation*}
\kappa_2(\mathcal{P})  \leq C^L \left( 1 + \frac{b_0 C_1 }{C-1} \right) (1+k_0 L) .
\end{equation*}
with $C=2(1+ a_0 b_0 C_1)$.
\begin{proof} 
Follows from Lemma \ref{Lem:multilevel_stability} and Theorem \ref{thm:k2bound}.
\end{proof}
\end{theorem}

\begin{remark} The upper bound in Lemma \ref{lem:alternative_stability} predicts an exponential increase
of the condition number with the number of levels $L+1$. Since we are
only interested in a moderate number of levels
$L+1=3$ or $4$, this may be acceptable. Moreover, our numerical results
below suggest that the bound is pessimistic. In our experiments, we observe 
$\kappa_2(\mathcal{P}) = O(L)$, which may actually be due to the lower bound.
\end{remark}

\subsection{Generalized Eigenproblems} \label{sec:GEVP}

The stability estimates in Definition~\ref{def:local_stability} are closely linked to properties
of the GEVP solved in each subdomain. This subsection establishes the necessary results.
In this section, $a : V \times V \to \mathbb{R}$ and $b : V \times V \to \mathbb{R}$ denote
generic symmetric and positive semi-definite bilinear forms on a $n$-dimensional vector space $V$. 
For a symmetric and positive semi-definite bilinear form $a$ denote by
$$\ker a = \{ v \in V \,:\, a(v,w)=0 \  \forall w \in V \}$$
the kernel of $a$. In the following, we make the important assumption
(to be verified below for the bilinear forms in \eqref{eq:theGEVP}) that the bilinear forms satisfy
\begin{equation}\label{eq:zero_kernel_intersection}
\ker a \cap \ker b = \{ 0 \} .
\end{equation}

\begin{definition}[Generalized Eigenvalue Problem]
We call $(\lambda,p) \in \mathbb{R} \cup \{ \infty \} \times V$ 
with $p \neq 0$, an eigenpair of $(a,b)$, if either $p \notin \ker b$ and
\begin{align}
\label{eq:GEVP}
a(p,v) = \lambda b(p,v) \qquad \forall v\in V
\end{align}
or $p \in \ker b$ and $\lambda = \infty$. For such a pair, $\lambda$ is called an eigenvalue 
and $p$ an eigenvector of $(a,b)$. The collection of all eigenvalues
(counted according to their 
geometric multiplicity) is called the spectrum of $(a,b)$.
\end{definition}

\begin{lemma}
The generalized eigenvalue problem \eqref{eq:GEVP} is non-defective, i.e. it has a full set of eigenvectors
with either $0\leq \lambda \in \mathbb{R}$ or $\lambda=+\infty$.
\begin{proof}
Assumption \eqref{eq:zero_kernel_intersection} implies that $a+b$ is a symmetric and positive definite bilinear form.
Employing a spectral transformation yields
\begin{align*}
a(p,v) = \lambda b(p,v) \quad
\Leftrightarrow \quad b(p,v) = \mu (a+b)(p,v) \quad \forall v\in V, \qquad \mu = \frac{1}{1+\lambda},
\end{align*}
with a symmetric and positive definite bilinear form on the right. Now standard spectral theory for this problem
shows that there exists a full set of eigenvectors with real nonnegative eigenvalues. It also
follows that $\mu\leq 1$. Now $\mu=0$ corresponds to $\lambda=+\infty$, $p\in \ker b$ and
$\mu=1$ corresponds to $\lambda=0$, $p\in \ker a$.
\end{proof}
\end{lemma}

From now on we assume that eigenpairs $(\lambda_k,p_k)$ are ordered by size, i.e. $\lambda_k \leq \lambda_{k+1}$.
With $r=\dim\ker a$ and $s=\dim\ker b$ the spectrum reads
\begin{equation*}
0 = \lambda_1 = \ldots = \lambda_r < \lambda_{r+1} \leq \ldots \leq \lambda_{n-s} < \lambda_{n-s+1} = \ldots
= \lambda_n = +\infty.
\end{equation*}
The eigenvectors $p_1,\ldots,p_{n-s}$, corresponding to the first $n-s$
eigenvalues, can be chosen to be simultaneously $b$-orthonormal and
$a$-orthogonal. The remaining eigenvectors $p_{n-s+1},\ldots,p_n \in\ker b$ are chosen to be $a$-orthogonal
and are also $b$-orthogonal, such that $\{p_1,\ldots,p_n\}$ forms a
basis of $V$.

The projection $\Pi_m : V \to V$ to
the first $m \leq n-s$ eigenvectors is defined by
\begin{equation} \label{eq:gevp_projection}
\Pi_m v := \sum_{k=1}^m b(v,p_k) p_k\,.
\end{equation}
Furthermore, we also introduce the induced semi-norms
\begin{align*}
|v|_a &= \sqrt{a(v,v)}, & |v|_b &= \sqrt{b(v,v)}.
\end{align*}

\begin{lemma}\label{lem:gevp_properties}
Let $a, b$ be positive semi-definite bilinear forms on $V$ with $\ker a \cap \ker b = \{ 0 \}$ and
$(\lambda_k,p_k)$, $k=1,\ldots,n$, eigenpairs of \eqref{eq:GEVP} with $p_k$, $k\leq n-s$, $b$-orthonormal
and $p_k$, $k>n-s$, $a$-orthogonal. Then the projection $\Pi_m$ satisfies the stability estimates
\begin{equation}\label{eq:a_stability}
|\Pi_m v|_a \leq |v|_a \qquad \text{and} \qquad |v - \Pi_m v|_a \leq |v|_a, \qquad \forall v \in V
\end{equation}
and for $m\geq r$
\begin{equation}\label{eq:ba_stability}
| v - \Pi_m v |_b^2 \leq \frac{1}{\lambda_{m+1}} | v - \Pi_m v |_a^2, \qquad \forall v \in V.
\end{equation}
\begin{proof}
The case $\ker b = \{0\}$ is proved in \cite[Lemma 2.11]{Geneo14}.
The proof of \eqref{eq:a_stability} in~\cite{Geneo14} extends without any
modification to the case $\ker b \neq \{0\}$, since 
the $p_k$ are still $a$-orthogonal and $m \leq n-s$.
For the proof of \eqref{eq:ba_stability} a small modification is
necessary. Let $v=\sum_{k=1}^n \alpha_k p_k$ be arbitrary. Then,
\begin{align*}
| v - \Pi_m v |_b^2 &= b\left( \sum_{k=m+1}^{n} \alpha_k p_k ,  \sum_{k=m+1}^{n} \alpha_k p_k \right)
= \sum_{k=m+1}^{n-s} \alpha_k^2 
= \sum_{k=m+1}^{n-s} \alpha_k^2 \frac{a(p_k,p_k)}{\lambda_k}\\
&\leq \frac{1}{\lambda_{m+1}} \sum_{k=m+1}^{n-s} \alpha_k^2 a(p_k,p_k) 
\leq \frac{1}{\lambda_{m+1}} \sum_{k=m+1}^{n} \alpha_k^2 a(p_k,p_k) \\
&= \frac{1}{\lambda_{m+1}} a\left( v - \Pi_m v  ,v - \Pi_m v \right) 
=  \frac{1}{\lambda_{m+1}} | v - \Pi_m v |_a^2 .
\end{align*}
\end{proof}
\end{lemma}
 
\begin{remark}
In \cite{Geneo14} the bilinear form $b$ in the GenEO method is positive semi-definite
but the authors did not include this case in their Lemma 2.11 but rather treat this fact in Lemmata
3.11, 3.16 and 3.18 for a special case. We think the assumption $\ker a \cap \ker b = \{ 0 \}$
is more natural and easy to prove in applications. In \cite[Lemma 2.3]{GenEO19} the authors
consider the GEVP for two symmetric positive semi-definite matrices $A$, $B$ with nontrivial
intersection $\ker A \cap \ker B \neq \{ 0\}$. This is not required in our analysis in the variational
setting. However, in the implementation a basis for the space $V =
V_l$ is required on each level.  It turns out, it is 
prohibitively expensive to construct such a basis on levels $l<L$ and only a generating system
is available. This then results in a GEVP with $\ker A \cap \ker B \neq \{ 0\}$. We refer to 
\cref{sec:Impl_GEVP} how to overcome this
problem. 
\end{remark}

\subsection{Application to Continuous Galerkin}

In this section, we consider the application to the solution of the scalar elliptic boundary value problem
\begin{subequations}\label{eq:ScalarBVP}
\begin{align}
-\nabla \cdot( K \nabla u) &= f &&\text{in $\Omega$},\\
u &= 0 &&\text{on $\Gamma_D\subseteq\partial\Omega$},\\
-(K \nabla u)\cdot \nu &= \psi &&\text{on $\Gamma_N=\partial\Omega\setminus\Gamma_D$} 
\end{align}
\end{subequations}
with Dirichlet and Neumann boundary conditions in a domain $\Omega\subset\mathbb{R}^d$.
The diffusion coefficient $K(x)\in\mathbb{R}^{d\times d}$ is symmetric and positive definite with eigenvalues bounded 
uniformly from above and below for all $x \in \Omega$.

We discretize \eqref{eq:ScalarBVP} with conforming finite elements on simplicial or hexahedral meshes \cite{Ern}
resulting in the weak formulation
\begin{equation*}
u_h \in V_h : \qquad a_h(u,v) = l(v) \qquad \forall v\in V_h
\end{equation*}
with
\begin{align*}
a_h(u,v) &= \sum_{\tau\in\mathcal{T}_h} a_\tau(u,v),  &
a_\tau(u,v) &= \int\limits_\tau (K\nabla u) \cdot \nabla v \,dx, &
l(v) = \int\limits_\Omega fv \,dx - \int\limits_{\Gamma_N} \psi v \,ds.
\end{align*}
The local bilinear forms $\overline a_{l,i}$ on the left side of the GEVP \eqref{eq:theGEVP} are then defined as 
\begin{equation}\label{eq:locala_cg}
\overline a_{l,i}(u,v) = \sum_{\tau\in\mathcal{T}_{h,l,i}} a_\tau(u,v).
\end{equation}

The following Lemma shows that the $\overline a_{l,i}$ define a symmetric positive semi-definite splitting
and the strengthened triangle inequality under the square holds.
\begin{lemma} 
Let the domain decomposition satisfy definition \ref{def:coloring} with $k_0$.
Then the local bilinear forms \eqref{eq:locala_cg}
satisfy definitions \ref{def:binomial} and \ref{def:posdefsplitting} with $a_0 = b_0 = k_0$.
\begin{proof}
Follows from any $\tau\in\mathcal{T}_h$ being in at most $k_0$ subdomains on any level $l$.
\end{proof}
\end{lemma}

Extending the results in \cite{Geneo14}, we show local stability
(Def. \ref{def:local_stability}) for three different right hand sides
in the GEVPs:
\begin{subequations} \label{eq:cg_rhs_blf}
\begin{align}
\overline b_{l,i}^\text{1}(u,v) &=  \overline a_{l,i}\left(\chi_{l,i} u , \chi_{l,i} v  \right) = a\left( e_{l,i}\chi_{l,i} u , e_{l,i}\chi_{l,i} v  \right), \\
\overline b_{l,i}^\text{2}(u,v) &=  \mathring{a}_{l,i}\left( \chi_{l,i} u , \chi_{l,i} v  \right), \label{eq:original_geneo}\\
\overline b_{l,i}^\text{3}(u,v) &=  \overline a_{l,i} \left( u -\chi_{l,i} u , v - \chi_{l,i} v \right),
\end{align}
\end{subequations}
where
\begin{equation*}
\mathring{a}_{l,i}(u,v)= \sum_{\tau\in\mathring{\mathcal{T}}_{h,l,i}} a_\tau(u,v) \quad\text{and}\quad
\mathring{\mathcal{T}}_{h,l,i} = \left\{ \tau \in \mathcal{T}_{h,l,i} : \tau\in\mathcal{T}_{h,l,j} \text{ for some } j\neq i\right\}.
\end{equation*}
The choice \eqref{eq:original_geneo} defines the original GenEO
method from \cite{Geneo14}. The other two choices are
easier to apply on the coarse levels $l<L$ and for discontinuous Galerkin (below).

\begin{lemma} The bilinear forms $\overline b_{l,i}^\alpha$, $\alpha=1,2,3$, satisfy $\ker \overline a_{l,i} \cap \ker \overline b_{l,i}^\alpha = \{0\}$.
\begin{proof} For subdomains touching the Dirichlet boundary, $\overline a_{l,i}$ is positive definite, so $\ker \overline a_{l,i} = \{0\}$.
For subdomains not touching the Dirichlet boundary, $\ker \overline
a_{l,i} = \mathspan\{ c\} $, where $c \equiv 1$  is the
constant one function. Now $\chi_{l,i}c$ and $(1-\chi_{l,i})c$ are not in $\ker \overline a_{l,i}$ 
Moreover, $\chi_{l,i}c$ is also not in  $\ker \mathring{a}_{l,i}$.
\end{proof}
\end{lemma}

Each choice of 
$\overline b_{l,i}^\alpha$, $\alpha=1,2,3$, gives rise to a projection operator
$\Pi^\alpha_{l,i,m}$, as defined in \eqref{eq:gevp_projection}, onto the
eigenvectors corresponding to the smallest $m$ 
eigenvalues in subdomain $i$ on level $l$.
Typically, for some threshold $\eta>0$, we will choose $m=m(\eta)$ such that $\lambda_{l,i,k}\leq \eta$ for all $k\leq m(\eta)$.
Making use of these projection operators, the stable two-level
splitting in Definition \ref{def:local_stability}, 
for a level $l$ function $v_l \in V_{h,l}$, can be achieved by choosing
\begin{align} \label{eq:levelwise_splitting_cg}
v_{l,i} &= \chi_{l,i} (I - \Pi^\alpha_{l,i,m(\eta)}) r_{l,i} v_l, &
v_{l-1} &= \sum_{i=1}^{P_l}  \chi_{l,i} \Pi^\alpha_{l,i,m(\eta)} r_{l,i} v_l .
\end{align}
Now we can state the main result of this subsection.

\begin{lemma} \label{lem:spectral1}
For $\alpha=1,2,3$, the splittings \eqref{eq:levelwise_splitting_cg}
using the projections $\Pi^\alpha_{l,i,m(\eta)}$ 
to the first $m(\eta)$ eigenvectors are locally stable with constants $C_1^\alpha$ given by
\begin{align*}
C_1^{1} &= \eta^{-1}, &
C_1^{2} &= 1+\eta^{-1}, &
C_1^{3} &= 2\left(1+\eta^{-1}\right).
\end{align*}
\begin{proof}
All three cases follow from Lemma \ref{lem:gevp_properties}. In the case $\alpha=1$, we obtain
\begin{align*}
\left\|v_{l,i}\right\|_{a_h}^2 &= \left\| \chi_{l,i} (I - \Pi^1_{l,i,m(\eta)}) r_{l,i} v_l \right\|_{a_h}^2
=  \left\| (I - \Pi^1_{l,i,m(\eta)}) r_{l,i} v_l\right\|_{\overline b^1_{l,i}}^2\\
&\leq \eta^{-1} \left\| (I - \Pi^1_{l,i,m(\eta)}) r_{l,i} v_l \right\|_{\overline a_{l,i}}^2
\leq  \eta^{-1} \left\| r_{l,i} v_l \right\|_{\overline a_{l,i}}^2.
\end{align*}
For $\alpha=2$, observe first that $(\chi_{l,i} v)|_\tau =
v|_\tau$\,, for all elements
$\tau\in\mathcal{T}_{h,l,i}\setminus\mathring{\mathcal{T}}_{h,l,i}$. Thus,
\begin{align*}
\bigl\|v_{l,i}&\bigr \|_{a_h}^2 = \sum_{\tau\in\mathring{\mathcal{T}}_{h,l,i}} a_\tau\left(v_{l,i},v_{l,i}\right)
+\sum_{\tau\in\mathcal{T}_{h,l,i}\setminus\mathring{\mathcal{T}}_{h,l,i}} a_\tau\left(v_{l,i},v_{l,i}\right)\\
&\leq \left| (I - \Pi^2_{l,i,m(\eta)}) r_{l,i} v_l
  \right|_{\overline b_{l,i}}^2
+ \left| (I - \Pi^2_{l,i,m(\eta)}) r_{l,i} v_{l}\right|_{\overline a_{l,i}}^2\\
&
\leq \Big( \eta^{-1} + 1\Big)  \left| (I - \Pi^2_{l,i,m(\eta)})
  r_{l,i} v_l   \right|_{\overline a_{l,i}}^2
\leq  \left( 1 + \eta^{-1} \right) \left| r_{l,i} v_l \right|_{\overline a_{l,i}}^2 .
\end{align*}
Finally, consider $\alpha=3$: Since $v_{l,i}$ is zero on
  $\partial \Omega_{l,i}$\,,
\begin{align*}
\left\|v_{l,i}\right\|_{a_h}^2 &= \left| \chi_{l,i} (I - \Pi^3_{l,i,m(\eta)}) r_{l,i} v_l  \right|_{\overline a_{l,i}}^2\\
&= \left| (I - \Pi^3_{l,i,m(\eta)}) r_{l,i} v_l - (I-\chi_{l,i}) (I - \Pi^3_{l,i,m(\eta)}) r_{l,i} v_l  \right|_{\overline a_{l,i}}^2
\end{align*}
\begin{align*}
&\leq 2\left| (I - \Pi^3_{l,i,m(\eta)}) r_{l,i} v_l\right|_{\overline a_{l,i}}^2 + 2\left|  (I - \Pi^3_{l,i,m(\eta)}) r_{l,i} v_l \right|_{\overline b_{l,i}}^2\\
&\leq 
2(1 +\eta^{-1} ) \left|  (I - \Pi^3_{l,i,m(\eta)}) r_{l,i} v_l \right|_{\overline a_{l,i}}^2
\leq  2(1+\eta^{-1}) \left| r_{l,i} v_l \right|_{\overline a_{l,i}}^2 .
\end{align*}
\end{proof}
\end{lemma}

\subsection{Application to Discontinuous Galerkin}

Here we consider the weighted symmetric interior penalty (WSIP) discontinuous Galerkin (DG)
method from \cite{Ern01042009} for the problem \eqref{eq:cg_rhs_blf}. 
The DG finite element space of degree $k$ on the mesh $\mathcal{T}_h$ is
\begin{equation}\label{eq:DG_fem_space}
V_h^\text{DG} = \{ v\in L_2(\Omega) : v|_\tau = p_\tau \circ \mu_{\tau}^{-1}, p_\tau \in\mathbb{P}_k\}
\end{equation}
where $\mathbb{P}_k$ is either the set of polynomials of total degree $k$ for simplices
or the set of polynomials of maximum degree $k$ for cuboid elements. 
A function $v\in V_h^\text{DG}$ is two-valued on an interior face $\gamma\in\mathcal{F}_h^I$ and
by $v^-$ we denote the restriction to $\gamma$ from $\tau^-(\gamma)$ and by $v^+$ the restriction from $\tau^+(\gamma)$.
For any point $x\in \gamma \in \mathcal{F}_h^I$ we define the jump and the weighted average
\begin{align*}
\llbracket v \rrbracket (x) &= v^-(x)-v^+(x), &
\{ v \}_\omega (x) &= \omega^- v^-(x) - \omega^+ v^+(x)
\end{align*}
for some weights $\omega^- + \omega^+ = 1$, $\omega^\pm \geq 0$.
A particular choice of the weights depending on the
absolute permeability tensor $K$ has been introduced in \cite{Ern01042009}. Assuming that
$K^\pm$ is constant on $\tau^\pm(\gamma)$, they set
$\omega^- = \delta_{K\nu}^+/(\delta_{K\nu}^- + \delta_{K\nu}^+)$ and
$\omega^+ = \delta_{K\nu}^-/(\delta_{K\nu}^- + \delta_{K\nu}^+)$ for
$\delta_{K\nu}^{\pm} =\nu_\gamma^T K^\pm \nu_\gamma$.
Finally, for any domain $Q \subset \Omega$ we set
\begin{equation*}
(v,w)_{0,Q} = \int_Q v w \ dx, \qquad \| v \|_{0,Q} = \sqrt{(v,v)_{0,Q}}.
\end{equation*}
The WSIP-DG method \cite{Ern01042009}
for numerically solving \eqref{eq:ScalarBVP} now reads
\begin{equation} \label{eq:DG_discrete_form}
u_h \in V_h^\text{DG} \quad:\quad a_h^\text{DG}(u_h,v) =  l_h^\text{DG}(v) \qquad \forall v \in V_h^\text{DG},
\end{equation}
where
\begin{equation} \label{eq:DG:blf}
a_h^\text{DG}(u,v) = 
\sum_{\tau\in\mathcal{T}_h} a_\tau(u,v) 
+ \sum_{\gamma\in\mathcal{F}_h^I}  a^I_\gamma(u,v)
+ \sum_{\gamma\in\mathcal{F}_h^D} a^D_\gamma(u,v)
\end{equation}
with $\mathcal{F}_h^D\subseteq \mathcal{F}_h^{\partial\Omega}$ denoting faces on $\Gamma_D$
and where 
\begin{align*}
a^I_\gamma(u,v) &= \sigma_\gamma(\llbracket u \rrbracket , \llbracket v \rrbracket  )_{0,\gamma}  
- ( \{K\nabla u\}_{\omega} \cdot \nu_\gamma , \llbracket  v \rrbracket   )_{0,\gamma}   
- ( \{K\nabla v\}_{\omega} \cdot \nu_\gamma , \llbracket  u \rrbracket   )_{0,\gamma}, \\
a^D_\gamma(u,v) &= \sigma_\gamma(u,v )_{0,\gamma} - ( (K\nabla u) \cdot \nu_\gamma , v  )_{0,\gamma}   
- ( (K\nabla v) \cdot \nu_\gamma , u  )_{0,\gamma} ,\\
l^\text{DG}(v) &= \sum_{\tau\in\mathcal{T}_h} (f,v)_{0,\tau} - \sum_{\gamma\in\mathcal{F}_h^N} (\psi,v)_{0,\gamma}
- \sum_{\gamma\in\mathcal{F}_h^D} \Bigl[  ( (K\nabla v) \cdot \nu_\gamma , g  )_{0,\gamma} 
- \sigma_\gamma(g,v )_{0,\gamma} \Bigr] .
\end{align*}
For each face (interior or boundary), $\sigma_\gamma>0$ defines a penalty parameter to be chosen and specified below.
The bilinear form $a_h^\text{DG}$ is symmetric and positive definite, provided the penalty parameters $\sigma_\gamma$ 
are chosen large enough \cite{Ern01042009}. 

We now prepare some results necessary for proving the requirements stated in Definitions  \ref{def:binomial} and \ref{def:posdefsplitting}.
\begin{lemma}\label{Lem:FluxEstimate}
Let $K(x)$ be a diffusion coefficient, constant on each element,
and choose $s\in\mathbb{R}$, $s>0$ an arbitrary number. Then the following estimates hold:
\begin{align*}
\Bigl| 2 ( (K\nabla v)\cdot\nu_\gamma, v)_{0,\gamma} \Bigr| &\leq
\frac{1}{s} a_{\tau^-(\gamma)}(v,v) +\theta_\gamma \| v \|_{0,\gamma}^2, &&\gamma\in\mathcal{F}_h^D,\\
\Bigl| 2 ( \{K\nabla v\}_\omega\cdot\nu_\gamma, \llbracket v \rrbracket )_{0,\gamma} \Bigr| &\leq
\frac{1}{s} (a_{\tau^-(\gamma)}(v,v) + a_{\tau^+(\gamma)}(v,v)) + \theta_\gamma \| \llbracket v \rrbracket \|_{0,\gamma}^2, &&\gamma\in\mathcal{F}_h^I,
\end{align*} 
with\vspace{-1ex}
\begin{equation*}
\theta_\gamma(s) = \begin{cases}
\displaystyle\frac{C_t^2 s \delta_{K\nu}}{2 h_\tau} & \gamma\in\mathcal{F}_h^D \\
\displaystyle\frac{C_t^2 s}{\min(h_{\tau^-},h_{\tau^+})} \frac{\delta_{K\nu}^- \delta_{K\nu}^+}{\delta_{K\nu}^- + \delta_{K\nu}^+} &\gamma\in\mathcal{F}_h^I
\end{cases}.
\end{equation*}
\begin{proof} Follows from the proof of  \cite[Lemma 3.1]{Ern01042009}. 
\end{proof}
\end{lemma}

The elementwise bilinear forms $a_\tau$ are positive semi-definite
and induce semi-norms $|v|_{a_\tau}=\sqrt{a_\tau(v,v)}$
that satisfy the triangle inequality $|v+w|_{a_\tau} \leq |v|_{a_\tau} +|w|_{a_\tau}$. A similar statement is not true
for the face bilinear forms $a_\gamma^{I}$ and $a_\gamma^{D}$, but one can prove the following estimate involving in
addition the elements adjacent to the face.
\begin{lemma} \label{Lem:triangle_alternative}
For any $s>0$ and for any penalty parameters $\sigma_\gamma > \theta_\gamma
= \theta_\gamma(s)$, let $v_1,\ldots,v_{k}$ be $k \ge 1$
arbitrary functions in $V_h^\text{DG}$. Then, the following bounds on
interior and boundary faces hold:
\begin{align*}
a^I_\gamma\Biggl(\sum_{i=1}^{\new{k}} v_i, \sum_{i=1}^{\new{k}} v_i\Biggr) 
&\leq {\new{k}} \frac{\sigma_\gamma + \new{\theta}_\gamma}{\sigma_\gamma-\new{\theta}_\gamma} \sum_{i=1}^{\new{k}}  a_\gamma^{I}(v_i,v_i)\\
&+ \frac{\new{k}}{\new{s}} \frac{2 \sigma_\gamma}{\sigma_\gamma-\new{\theta}_\gamma} \sum_{i=1}^{\new{k}} a_{\tau^-(\gamma)}\left( v_i, v_i\right)
+ \frac{\new{k}}{\new{s}} \frac{2 \sigma_\gamma}{\sigma_\gamma-\new{\theta}_\gamma} \sum_{i=1}^{\new{k}} a_{\tau^+(\gamma)}\left( v_i, v_i\right),\\
a^D_\gamma\left(\sum_{i=1}^{\new{k}} v_i, \sum_{i=1}^{\new{k}} v_i\right) 
&\leq {\new{k}} \frac{\sigma_\gamma + \new{\theta}_\gamma}{\sigma_\gamma-\new{\theta}_\gamma} \sum_{i=1}^{\new{k}}  a_\gamma^{D}(v_i,v_i)
+ \frac{\new{k}}{\new{s}} \frac{2 \sigma_\gamma}{\sigma_\gamma-\new{\theta}_\gamma} \sum_{i=1}^{\new{k}} a_{\tau^-(\gamma)}\left( v_i, v_i\right) .
\end{align*} 
\begin{proof} Using Lemma \ref{Lem:FluxEstimate},
we obtain, for a single function $v$ and an interior face~$\gamma$
with adjacent elements $\tau^-=\tau^-(\gamma)$ and
$\tau^+=\tau^+(\gamma)$, that
\begin{align*}
a^I_\gamma(v,v) &= \sigma_\gamma \| \llbracket v \rrbracket\|_{0,\gamma}^2 - 2 ( \{K\nabla v\}_\omega\cdot\nu_\gamma, \llbracket v \rrbracket )_{0,\gamma}\\
&\leq (\sigma_\gamma + \new{\theta}_\gamma) \| \llbracket v \rrbracket\|_{0,\gamma}^2 
+  \frac{1}{\new{s}} a_{\tau^-}(v,v) + \frac{1}{\new{s}} a_{\tau^+}(v,v)
\end{align*}
as well as
\begin{align*}
a^I_\gamma(v,v) &= \sigma_\gamma \| \llbracket v \rrbracket\|_{0,\gamma}^2 - 2 ( \{K\nabla v\}_\omega\cdot\nu_\gamma, \llbracket v \rrbracket )_{0,\gamma}\\
&\geq (\sigma_\gamma-\new{\theta}_\gamma) \| \llbracket v \rrbracket\|_{0,\gamma}^2 -
\frac{1}{\new{s}} a_{\tau^-}(v,v) - \frac{1}{\new{s}} a_{\tau^+}(v,v) \\
\Leftrightarrow\quad  (\sigma_\gamma-\new{\theta}_\gamma)  \| \llbracket v \rrbracket\|_{0,\gamma}^2 
&\leq a^I_\gamma(v,v) + \frac{1}{\new{s}} a_{\tau^-}(v,v) + \frac{1}{\new{s}} a_{\tau^+}(v,v) .
\end{align*}
Now observe that $\| \llbracket v \rrbracket\|_{0,\gamma}^2$ and $a_\tau(v,v)$ are (semi-)norms for which the
triangle inquality holds and $|\sum_{i=1}^{\new{k}} v_i |^2 \leq
{\new{k}}\sum_{i=1}^{\new{k}} | v_i |^2$. Then,
\begin{align*}
a^I_\gamma&\Biggl(\sum_{i=1}^{\new{k}} v_i,  \sum_{i=1}^{\new{k}} v_i\Biggr)
\leq (\sigma_\gamma + \new{\theta}_\gamma) \left\| \left\llbracket \sum_{i=1}^{\new{k}} v_i \right\rrbracket\right\|_{0,\gamma}^2 \\
&\qquad\qquad+  \frac{1}{\new{s}} a_{\tau^-}\left(\sum_{i=1}^{\new{k}} v_i, \sum_{i=1}^{\new{k}} v_i\right) 
+ \frac{1}{\new{s}} a_{\tau^+}\left(\sum_{i=1}^{\new{k}} v_i, \sum_{i=1}^{\new{k}} v_i\right) \\
&\leq (\sigma_\gamma + \new{\theta}_\gamma) {\new{k}} \sum_{i=1}^{\new{k}}  \| \left\llbracket v_i \right\rrbracket\|_{0,\gamma}^2 
 + \frac{{\new{k}}}{\new{s}} \sum_{i=1}^{\new{k}} a_{\tau^-}( v_i, v_i) 
+ \frac{{\new{k}}}{\new{s}} \sum_{i=1}^{\new{k}} a_{\tau^+}\left(v_i, v_i\right) 
\end{align*}
\begin{align*}
\qquad &\leq {\new{k}} \frac{\sigma_\gamma + \new{\theta}_\gamma}{\sigma_\gamma-\new{\theta}_\gamma}  
\Biggl( \sum_{i=1}^{\new{k}}  a^I_\gamma(v_i,v_i) + \frac{1}{\new{s}} \sum_{i=1}^{\new{k}} a_{\tau^-}(v_i,v_i) 
+ \frac{1}{\new{s}} \sum_{i=1}^{\new{k}} a_{\tau^+}(v_i,v_i)  \Biggr)\\
&\qquad\qquad + \frac{{\new{k}}}{\new{s}} \sum_{i=1}^{\new{k}} a_{\tau^-}( v_i, v_i) 
+ \frac{{\new{k}}}{\new{s}} \sum_{i=1}^{\new{k}} a_{\tau^+}(v_i, v_i) 
\end{align*}
from which the result is obtained.
Boundary faces are treated similarly.
\end{proof}
\end{lemma}

For the DG method the left-hand side bilinear form in the GEVP is defined as
\begin{equation} \label{eq:locala_dg} 
\overline a_{h,l,i}^\text{DG}(u,v) = 
\sum_{\tau\in\mathcal{T}_{h,l,i}} a_\tau(u,v) 
+ \sum_{\gamma\in\mathcal{F}_{h,l,i}^{I}}  a^I_\gamma(u,v)
+ \sum_{\gamma\in\mathcal{F}_{h,l,i}^{D}} a^D_\gamma(u,v),
\end{equation}
i.e. \new{only faces that are interior to $\Omega_{l,i}$ are used,
  which corresponds to} Neumann boundary conditions on  $\partial\Omega_{l,i}\cap\Omega$.

\begin{lemma} \label{lem:summability_DG}
Let the domain decomposition satisfy Definition \ref{def:coloring}
with $k_0 \in \mathbb{N}$.
Let $m_{\mathcal F}$ be the maximum number of faces of an element, choose
$s>0$ and 
$\sigma_\gamma>\theta_\gamma(s)$.
Then the local bilinear forms \eqref{eq:locala_dg}
satisfy Definitions \ref{def:binomial} and~\ref{def:posdefsplitting} with $b_0 = k_0$ and
\begin{align*}
a_0 &=  k_0 \max\left( 1 + \frac{m_{\mathcal F}}{s} \max_{\gamma\in\mathcal{F}_h^I\cup \gamma\in\mathcal{F}_h^D} \frac{2 \sigma_\gamma}{\sigma_\gamma-\theta_\gamma},
\max_{\gamma\in\mathcal{F}_h^I\cup \gamma\in\mathcal{F}_h^D}\frac{\sigma_\gamma + \theta_\gamma}{\sigma_\gamma-\theta_\gamma} \right) .
\end{align*}
\begin{proof}
The proof is based on  Lemma \ref{Lem:triangle_alternative}, choosing
$k = k_0$. Set $J_{l,\tau} = \{ i : \tau \in \mathcal{T}_{h,l,i}\}$,
 $J_{l,\gamma} = J_{l,\tau^-(\gamma)} \cup J_{l,\tau^+(\gamma)}$ for interior faces
 and $J_{l,\gamma} = J_{l,\tau^-(\gamma)}$ for boundary faces. Observe that $|J_{l,\tau}|\leq k_0$ and $|J_{l,\gamma}|\leq k_0$
and estimate
\begin{align*}
\Biggl\| & \sum_{i=1}^{P_l}  v_{l,i} \Biggr\|_{a_h^\text{DG}}^2 = a_h^\text{DG}\left(\sum_{i=1}^{P_l} v_{l,i}, \sum_{i=1}^{P_l} v_{l,i}\right) 
\ = \ \sum_{\tau\in\mathcal{T}_h} a_\tau\left(\sum_{i\in J_{l,\tau}} v_{l,i},\sum_{i\in J_{l,\tau}} v_{l,i}\right) \\
&\qquad\qquad + \sum_{\gamma\in\mathcal{F}_h^I}  a^I_\gamma\left(\sum_{i\in J_{l,\gamma}} v_{l,i},\sum_{i\in J_{l,\gamma}} v_{l,i}\right)
+ \sum_{\gamma\in\mathcal{F}_h^D} a^D_\gamma\left(\sum_{i\in J_{l,\gamma}} v_{l,i},\sum_{i\in J_{l,\gamma}} v_{l,i}\right) \\ 
&\leq k_0 \sum_{\tau\in\mathcal{T}_h} \sum_{i\in J_{l,\tau}} a_\tau\left( v_{l,i}, v_{l,i}\right) 
+ \sum_{\gamma\in\mathcal{F}_h^I} 
\Biggl[ k_0 \frac{\sigma_\gamma + \theta_\gamma}{\sigma_\gamma-\new{\theta}_\gamma} \sum_{i\in J_{l,\gamma}}  a_\gamma(v_{l,i},v_{l,i})
\end{align*}\vspace{-5ex}
\begin{align*}
&\qquad + \frac{k_0}{\new{s}} \frac{2 \sigma_\gamma}{\sigma_\gamma-\new{\theta}_\gamma} \sum_{i\in J_{l,\gamma}} a_{\tau^-(\gamma)}\left( v_{l,i}, v_{l,i}\right)
+ \frac{k_0}{\new{s}} \frac{2 \sigma_\gamma}{\sigma_\gamma-\new{\theta}_\gamma} \sum_{i\in J_{l,\gamma}} a_{\tau^+(\gamma)}\left( v_{l,i}, v_{l,i}\right)\Biggr] \\
&\qquad + \sum_{\gamma\in\mathcal{F}_h^D} 
\Biggl[ k_0 \frac{\sigma_\gamma + \new{\theta}_\gamma}{\sigma_\gamma-\new{\theta}_\gamma} \sum_{i\in J_{l,\gamma}}  a_\gamma(v_{l,i},v_{l,i})
+ \frac{k_0}{\new{s}} \frac{2 \sigma_\gamma}{\sigma_\gamma-\new{\theta}_\gamma} \sum_{i\in J_{l,\gamma}} a_{\tau^-(\gamma)}\left( v_{l,i}, v_{l,i}\right)\Biggr]
\end{align*}
\begin{align*}
\quad&\leq a_0 \left( \sum_{\tau\in\mathcal{T}_h} \sum_{i\in J_{l,\tau}} a_\tau\left( v_{l,i}, v_{l,i}\right) 
+ \sum_{\gamma\in\mathcal{F}_h^I} \sum_{i\in J_{l,\gamma}}  a_\gamma(v_{l,i},v_{l,i})
+ \sum_{\gamma\in\mathcal{F}_h^D} \sum_{i\in J_{l,\gamma}}  a_\gamma(v_{l,i},v_{l,i}) \right)\\
&= a_0 \left( \sum_{\tau\in\mathcal{T}_h} \sum_{i=1}^{P_l} a_\tau\left( v_{l,i}, v_{l,i}\right) 
+ \sum_{\gamma\in\mathcal{F}_h^I} \sum_{i=1}^{P_l}  a_\gamma(v_{l,i},v_{l,i})
+ \sum_{\gamma\in\mathcal{F}_h^D} \sum_{i=1}^{P_l}  a_\gamma(v_{l,i},v_{l,i}) \right)\\
&= a_0 \sum_{i=1}^{P_l} a\left(v_{l,i},v_{l,i}\right) = a_0 \sum_{i=1}^{P_l} \left\| v_{l,i} \right\|_{a_h^\text{DG}}^2 . 
\end{align*}

Now consider the SPSD splitting property from Definition \ref{def:posdefsplitting}.
For any interior face $\gamma\in\mathcal{F}_{h}^I$, 
set $J_{l,\gamma}^\ast = J_{l,\tau^-(\gamma)} \cap J_{l,\tau^+(\gamma)} \subset J_{l,\gamma}$
and observe that $J_{l,\gamma}^\ast\neq\emptyset$ due to the second condition in Definition \ref{def:coloring}. Then
\begin{align*}
\sum_{i=1}^{P_l} \left| r_{l,i} v_l \right|_{\overline a_i}^2 
&= \sum_{i=1}^{{P_l}} \Biggl( \sum_{\tau\in\mathcal{T}_{h,l,i}} a_\tau(v_l,v_l) 
+ \sum_{\gamma\in\mathcal{F}_{h,l,i}^{I}}  a^I_\gamma(v_l,v_l)
+ \sum_{\gamma\in\mathcal{F}_{h,l,i}^{D}} a^D_\gamma(v_l,v_l) \Biggr)\\
&=  \sum_{\tau\in\mathcal{T}_{h}}\sum_{i\in J_\tau} a_\tau(v_l,v_l) 
+ \sum_{\gamma\in\mathcal{F}_{h}^I} \sum_{i\in J^\ast_\gamma}  a^I_\gamma(v_l,v_l)
+ \sum_{\gamma\in\mathcal{F}_{h}^D} \sum_{i\in J_\gamma} a^D_\gamma(v_l,v_l)\\
&\leq k_0 \left( \sum_{\tau\in\mathcal{T}_{h}} a_\tau(v_l,v_l) 
+ \sum_{\gamma\in\mathcal{F}_{h}^I}  a^I_\gamma(v_l,v_l)
+ \sum_{\gamma\in\mathcal{F}_{h}^D} a^D_\gamma(v_l,v_l)\right)
= k_0 \left\| v_l \right\|_{a_h^\text{DG}}^2 .
\end{align*}
\end{proof}
\end{lemma}

As in the continuous Galerkin case, the local stability in Definition \ref{def:local_stability} 
is ensured by solving appropriate GEVPs. Three possible right-hand side bilinear forms are
\begin{subequations} \label{eq:dg_rhs_blf}
\begin{align}
\overline b_{l,i}^\text{DG,1}(u,v) &= \overline a^\text{DG}_{h,l,i}\left(\chi_{l,i} u ,\chi_{l,i} v  \right) = a_h^\text{DG}\left( e_{l,i}\chi_{l,i} u , e_{l,i}\chi_{l,i} v  \right), \\
\overline b_{l,i}^\text{DG,2}(u,v) &=  \mathring{a}_{h,l,i}^\text{DG} \left( \chi_{l,i} u , \chi_{l,i} v  \right),\\
\overline b_{l,i}^\text{DG,3}(u,v) &=  \overline a^\text{DG}_{h,l,i} \left( u -\chi_{l,i} u , v - \chi_{l,i} v \right).
\end{align}
\end{subequations}
resulting in corresponding projection operators $\Pi^\text{DG,
  $\alpha$}_{l,i,m(\eta)}$. 

For $\alpha=2$, the
original GenEO method, the bilinear form is given by
\begin{equation*}
\mathring{a}_{h,l,i}^\text{DG}(u,v) = 
\sum_{\tau\in\mathring{\mathcal{T}}_{h,l,i}} a_\tau(u,v) 
+ \sum_{\gamma\in\mathring{\mathcal{F}}_{h,l,i}^{I}}  a^I_\gamma(u,v)
+ \sum_{\gamma\in\mathcal{F}_{h,l,i}^{D}} a^D_\gamma(u,v)
\end{equation*}
with $\mathring{\mathcal{F}}_{h,l,i}^{I}=\{\gamma\in\mathcal{F}_{h,l,i}^{I} : 
\tau^-(\gamma)\in\subset\mathring{\mathcal{T}}_{h,l,i} \vee \tau^+(\gamma)\in\subset\mathring{\mathcal{T}}_{h,l,i} \}$

\begin{lemma} \label{lem:spectral_dg}
For $\alpha=1,2,3$, using 
the projection operators $\Pi^\text{DG, $\alpha$}_{l,i,m(\eta)}$ 
to the first $m(\eta)$ eigenvectors, the splittings \eqref{eq:levelwise_splitting_cg} satisfy Definition \ref{def:local_stability} with 
\begin{align*}
C_1^\text{DG,1} &= \eta^{-1}, &
C_1^\text{DG,2} &= C_2+\eta^{-1}, &
C_1^\text{DG,3} &= 2\left(1+\eta^{-1}\right).
\end{align*}
where 
\begin{align*}
C_2& =\max\{1+C_\tau,C_\gamma^I,C_\gamma^D\}, & \qquad C_\tau & = \max_{\tau\in\mathcal{T}_{h,l,i}}
                                              \frac{(\sigma_\gamma+\new{\theta}_\gamma)\new{m_{\mathcal
                                              F}}}{(\sigma_\gamma-\new{\theta}_\gamma)
                                              \new{s}},\\
C_\gamma^I & = \max_{\gamma\in\mathcal{F}_{h,l,i}^I}
  \max(1,\frac{\sigma_\gamma+\new{\theta}_\gamma}{\sigma_\gamma-\new{\theta}_\gamma}), & 
\qquad C_\gamma^D & = \max_{\gamma\in\mathcal{F}_{h,l,i}^D}
  \frac{\sigma_\gamma+\new{\theta}_\gamma}{\sigma_\gamma-\new{\theta}_\gamma}.
\end{align*}
\begin{proof} For $\alpha=1,3$ the proof is the same as in Lemma \ref{lem:spectral1}.
For $\alpha=2$, we abbreviate $w_{l,i} = (I - \Pi^\text{DG,2}_{l,i,m(\eta)}) r_{l,i} v_{l}$
and $v_{l,i} = \chi_{l,i}w_{l,i}$ and observe as before
\begin{align*}
\sum_{\tau\in\mathcal{T}_{h,l,i}\setminus\mathring{\mathcal{T}}_{h,l,i}} a_\tau(v_{l,i},v_{l,i}) 
&= \sum_{\tau\in\mathcal{T}_{h,l,i}\setminus\mathring{\mathcal{T}}_{h,l,i}} a_\tau(w_{l,i},w_{l,i})
\leq \sum_{\tau\in\mathcal{T}_{h,l,i}} a_\tau(w_{l,i},w_{l,i}) .
\end{align*}
For the skeleton terms, we observe
\begin{align*}
&\new{\sum_{\gamma\in\mathcal{F}_{h,l,i}^{I}\setminus\mathring{\mathcal{F}}_{h,l,i}^{I}}  a^I_\gamma(v_{l,i},v_{l,i}) \le
\sum_{\gamma\in\mathcal{F}_{h,l,i}^{I}\setminus\mathring{\mathcal{F}}_{h,l,i}^{I}}  a^I_\gamma(w_{l,i},w_{l,i})}\\
& 
\qquad\qquad \new{+ \sum_{\gamma\in\mathring{\mathcal{F}}_{h,l,i}^{I}}  (\sigma_\gamma+\new{\theta}_\gamma) \|\llbracket w_{l,i} \rrbracket\|^2_{0,\gamma}
+ \sum_{\gamma\in\mathcal{F}_{h,l,i}^{D}}  (\sigma_\gamma+\new{\theta}_\gamma) \| w_{l,i} \|^2_{0,\gamma}}\\
&\qquad\leq \sum_{\gamma\in\mathcal{F}_{h,l,i}^{I}\setminus\mathring{\mathcal{F}}_{h,l,i}^{I}}  a^I_\gamma(w_{l,i},w_{l,i})
+  \sum_{\gamma\in\mathring{\mathcal{F}}_{h,l,i}^{I}}  \frac{\sigma_\gamma+\new{\theta}_\gamma}{\sigma_\gamma-\new{\theta}_\gamma} a^I_\gamma(w_{l,i},w_{l,i})\\
&\qquad\qquad + \sum_{\gamma\in\mathring{\mathcal{F}}_{h,l,i}^{I}} \frac{\sigma_\gamma+\new{\theta}_\gamma}{(\sigma_\gamma-\new{\theta}_\gamma) \new{s}}(a_{\tau^-}(w_{l,i},w_{l,i}) + a_{\tau^+}(w_{l,i},w_{l,i}))\\
&\qquad\qquad + \sum_{\gamma\in\mathcal{F}_{h,l,i}^{D}} \frac{\sigma_\gamma+\new{\theta}_\gamma}{\sigma_\gamma-\new{\theta}_\gamma} a^D_\gamma(w_{l,i},w_{l,i})
+ \sum_{\gamma\in\mathcal{F}_{h,l,i}^{D}}  \frac{\sigma_\gamma+\new{\theta}_\gamma}{(\sigma_\gamma-\new{\theta}_\gamma) \new{s}}a_{\tau^-}(w_{l,i},w_{l,i}) \\
&\qquad\leq C_\tau \sum_{\tau\in\mathcal{T}_{h,l,i}} a_\tau(w_{l,i},w_{l,i}) 
+ C_\gamma^D \sum_{\gamma\in\mathcal{F}_{h,l,i}^{D}}  a^D_\gamma(w_{l,i},w_{l,i})
+ C_\gamma^I \sum_{\gamma\in\mathcal{F}_{h,l,i}^{I}}  a^I_\gamma(w_{l,i},w_{l,i})
\end{align*}
Using these two results we estimate
\begin{align*}
\| v_{l,i} \|^2_{a_h^\text{DG}} &= 
\sum_{\tau\in\mathcal{T}_{h,l,i}} a_\tau(v_{l,i},v_{l,i}) 
+ \sum_{\gamma\in\mathcal{F}_{h,l,i}^{I}}  a^I_\gamma(v_{l,i},v_{l,i})
+ \sum_{\gamma\in\mathcal{F}_{h,l,i}^{D}} a^D_\gamma(v_{l,i},v_{l,i})\\
&= \mathring{a}_{h,l,i}^\text{DG}(v_{l,i},v_{l,i}) + 
\sum_{\tau\in\mathcal{T}_{h,l,i}\setminus\mathring{\mathcal{T}}_{h,l,i}} a_\tau(v_{l,i},v_{l,i}) 
+ \sum_{\gamma\in\mathcal{F}_{h,l,i}^{I}\setminus\mathring{\mathcal{F}}_{h,l,i}^{I}}  a^I_\gamma(v_{l,i},v_{l,i}) \\
&\leq \mathring{a}_{h,l,i}^\text{DG}(v_{l,i},v_{l,i}) 
+ C_2 \overline a_{h,l,i}^\text{DG}(w_{l,i},w_{l,i}) \\
&= \left| (I - \Pi^\text{DG,2}_{l,i,m(\eta)}) r_{l,i} v_l   \right|^2_{\overline b_{l,i}^\text{DG,2}} + C_2 \left| (I - \Pi^\text{DG,2}_{l,i,m(\eta)}) r_{l,i} v_{l}\right|^2_{\overline a_{h,l,i}^\text{DG}}\\
&\leq \eta^{-1} \left| (I - \Pi^2_{l,i,m(\eta)}) r_{l,i} v_l   \right|_{\overline a_{h,l,i}^\text{DG}}^2 + C_2\left| r_{l,i} v_{l}\right|_{\overline a_{h,l,i}^\text{DG}}^2\\
&\leq  \left( C_2 + \eta^{-1} \right) \left| r_{l,i} v_l \right|_{\overline a_{h,l,i}^\text{DG}}^2 .
\end{align*}
\end{proof}
\end{lemma}

\section{Implementation}
\label{sec:impl}

The multilevel spectral domain decomposition preconditioner described in this paper has been implemented
within the DUNE software framework\footnote{www.dune-project.org} \cite{Dune1,Duneneu}
in a sequential setting. Parallel runtimes reported e.g. in Table \ref{fig:islands3d} below are estimated from
sequential runs by taking the maximum over the times needed for the computations in each subdomain.

\subsection{Patch-wise stiffness matrices}

The implementation is fully algebraic in the sense that only input on the finest level is required. This input is
in the form of stiffness matrices assembled on certain nonoverlapping
sets of elements, which we call patch matrices.
Recall that by $J_{l,\tau} = \{i : \tau \in
\mathcal{T}_{h,l,i}\}$ we denoted the set of subdomain numbers that contain element
$\tau$ on level $l$. Each $\sigma\subset\{1,\ldots,P_l\}$ gives rise to a \textit{volume patch}
$\mathcal{T}_{l,\sigma} = \{\tau\in\mathcal{T} : J_{l,\tau} =
\sigma\}$. The volume patch matrix $A_{l,\sigma}$ contains all
contributions from elements in $\mathcal{T}_{l,\sigma}$.
In DG methods, the face terms need to be considered in addition. Boundary face contributions are assembled to the 
volume patch matrix of the corresponding element. Interior face contributions are assembled to a volume patch matrix if both 
adjacent elements belong to the same patch. 
Only if the two elements adjacent to a face belong to two different patches,
the contribution of that face is assembled to a seperate \textit{skeleton patch matrix} $A_{l,\sigma^-,\sigma^+}$,
collecting all contributions from the faces $\mathcal{F}_{l,\sigma^-,\sigma^+} = \{ \gamma\in\mathcal{F}_h^I : 
J_{l,\tau^-(\gamma)}=\sigma^- \wedge J_{l,\tau^+(\gamma)}=\sigma^+\}$.
The preconditioner gets volume patch matrices and skeleton patch matrices as input.
From this information all the relevant subdomain matrices on all levels can be computed.

\subsection{Solving the GEVPs}\label{sec:Impl_GEVP}

Implementing the multilevel spectral domain decomposition method requires a basis representation
for the spaces $\overline V_{h,l,i}$ introduced in
\eqref{eq:def_neumann_space}, which are used in the GEVP \eqref{eq:theGEVP}.
Consider a subdomain $i$ on level $l<L$ which is made up of the subdomains $j\in J_{l,i}$ from level $l+1$.
$\overline V_{h,l,i}$ is constructed from eigenfunctions computed in all the subdomains 
$j \in \overline J_{l,i} = \{ j^\prime \in\{1,\ldots,P_{l+1}\} : \mathcal{T}_{h,l+1,j^\prime} \cap \mathcal{T}_{h,l,i} \neq \emptyset \}$
in the following way:
\begin{equation*}
\overline V_{h,l,i} = \mathspan\{ \phi_{l,j,k} : j\in J_{l,i}, \lambda_{l,j,k}<\eta \}
\oplus \mathspan\{ r_{l,i} \phi_{l,j,k} : j\in \overline J_{l,i} \setminus J_{l,i}, \lambda_{l,j,k}<\eta \} .
\end{equation*}
While the functions in the first set have support in $\Omega_{l,i}$ and are linearly independent, 
the functions in the second set are restrictions to $\Omega_{l,i}$ of basis
  functions from neighbouring subdomains, which typically only form a
generating system and are linearly dependent. 
A 
basis is not cheaply available.
When the GEVP \eqref{eq:theGEVP}
is assembled with this generating system, it leads to algebraic eigenvalue problems
\begin{equation}
A_{l,i} x_{l,i,k} = \lambda_{l,i,k} B_{l,i} x_{l,i,k}
\end{equation}
with $\ker A_{l,i} \cap \ker B_{l,i} \neq \{0\}$. This is
  outside the scope of the theory presented in \cref{sec:GEVP} and we
overcome this problem as follows. Let $A_{l,i} = L_{l,i} D_{l,i} L^T_{l,i} $ be an $LDL^T$ factorization of
$A_{l,i}$ and let $D_{l,i,\epsilon}$ be a regularized version of the diagonal matrix $D_{l,i}$ where zeroes on the diagonal
are replaced by $0<\epsilon\ll 1$. Then, consider the spectral transformation 
\begin{equation}\label{eq:evp_reg}
L_{l,i}^{-T} D_{l,i,\epsilon}^{-1} L_{l,i}^{-1} B_{l,i}  x_{l,i,k} = \mu_{l,i,k} x_{l,i,k}
\end{equation}
with $\mu_{l,i,k} =  \lambda_{l,j,k}^{-1}$, where we are now
interested in the \textit{largest} eigenvalues of \eqref{eq:evp_reg}.
Crucially, all $x\in \ker B_{l,i}$ are eigenvectors
corresponding to $\mu=0$, and this includes obviously $\ker
  A_{l,i} \cap \ker B_{l,i}$. On the other hand, vectors $x\in\ker A_{l,i} \cap \range  B_{l,i}$
'pass' the matrix $B_{l,i}$ and lead to very large eigenvalues of order $\mu=\epsilon^{-1}$.
Important for this method to work in practice is that vectors $x\in \ker B_{l,i}$ give $B_{l,i} x = 0$
also in finite precision. 

\subsection{Software and hardware used}

For solving the GEVP, we use Arpack~\cite{ArpackUG}
through the Arpack++ wrapper
in symmetric shift-invert mode. As subdomain solver, we use Cholmod \cite{Cholmod2008}
in the iteration phase and UMFPack \cite{UMFPack} in the eigenvalue solver. As graph partitioner,
we use ParMetis \cite{KARYPIS199896}. 
Run-times reported below are in seconds and were obtained on
an Intel(R) Xeon(R) Silver 4114 CPU @ 2.20GHz.

\section{Numerical Results}
\label{sec:results}

We test the new multilevel spectral domain decomposition
  preconditioners within a Krylov iteration (Conjugate Gradients or GMRES) to solve \eqref{eq:linear_system}.
In all examples, we stop the computation when $\|b-Ax^m\| < 10^{-8} \| b- Ax^0\|$ and report 
the number of iterations $\#IT=m$ needed.

\subsection{Islands Problem}

\begin{figure}[tbp]
\begin{center}
\includegraphics[width=0.4\textwidth]{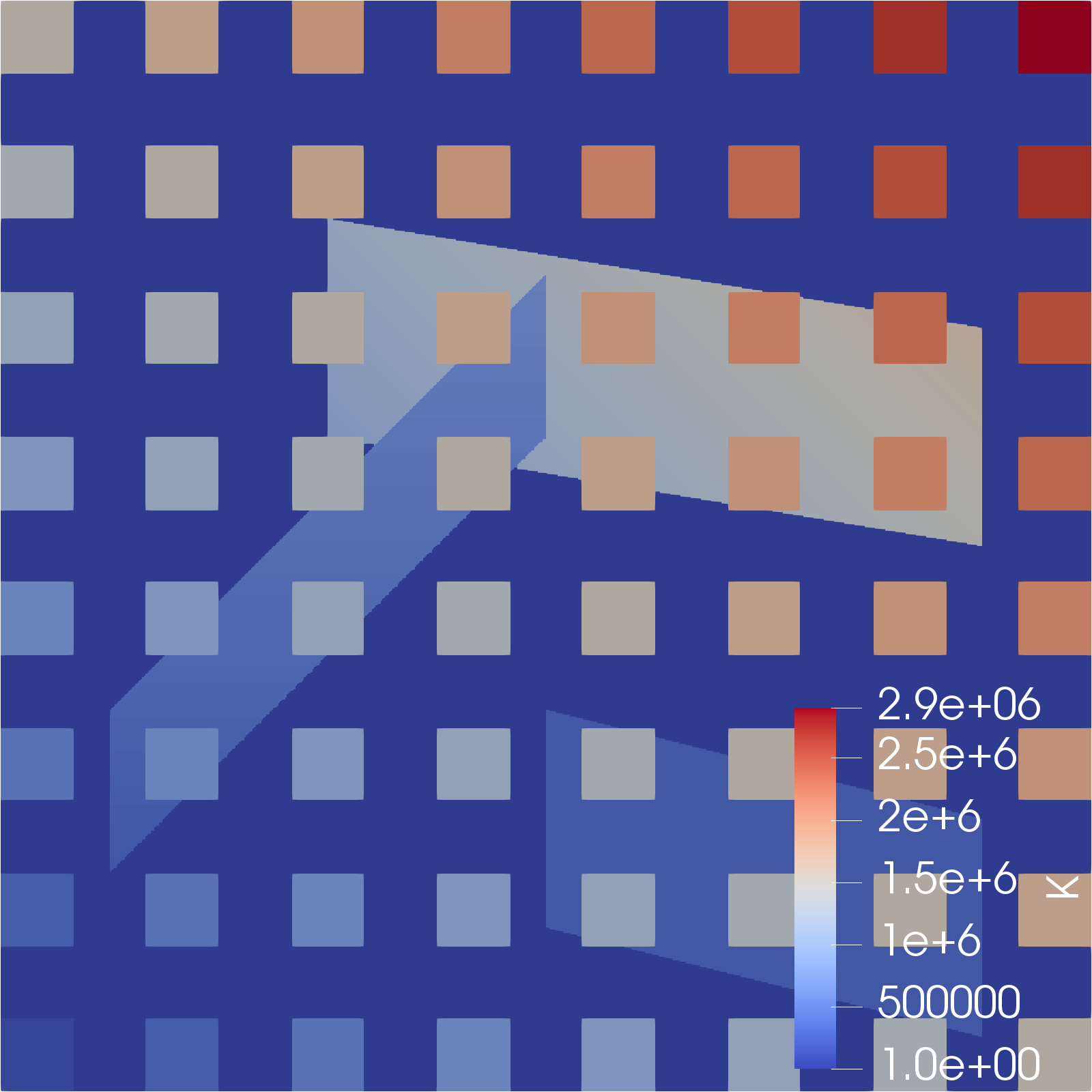}\hspace{3mm}
\includegraphics[width=0.4\textwidth]{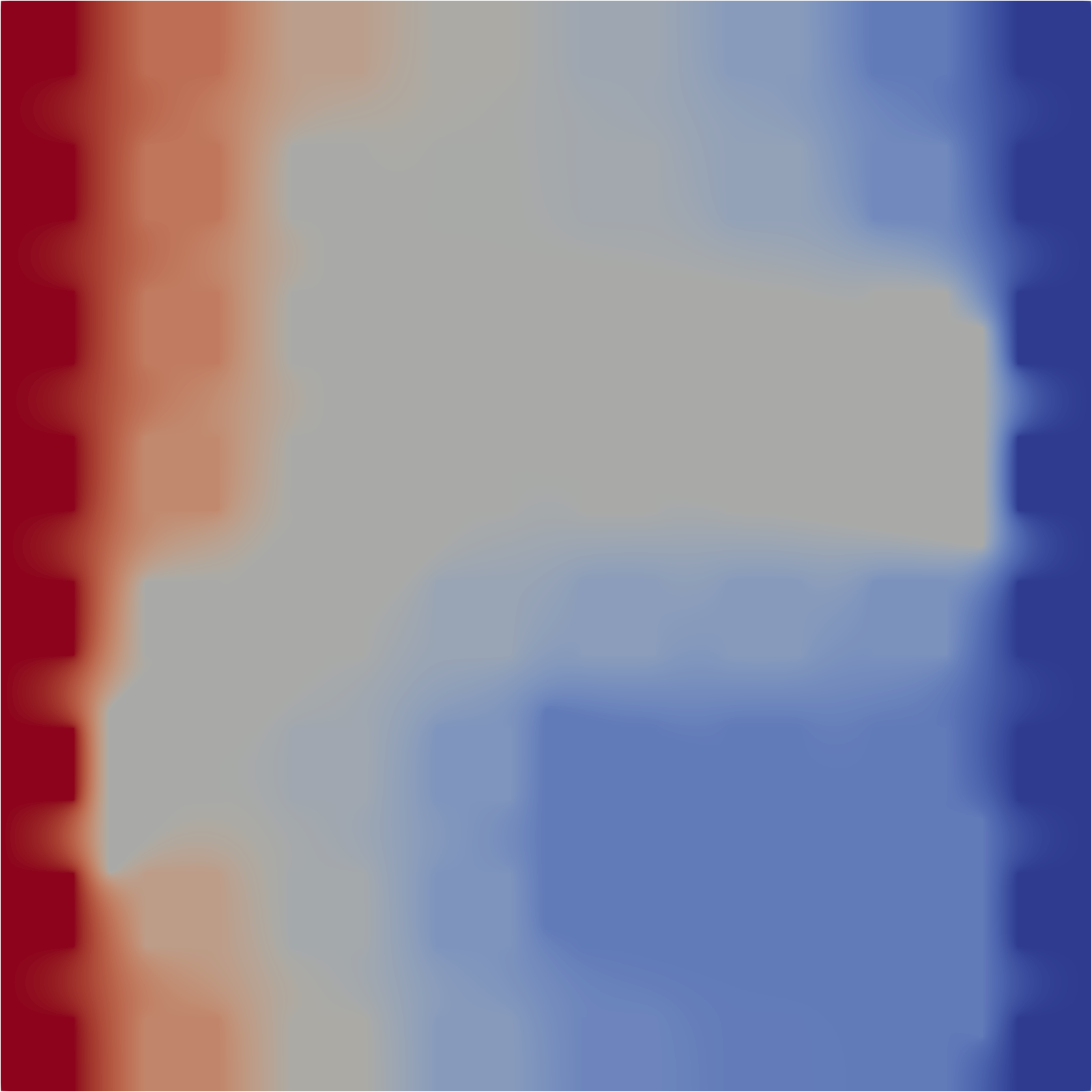}
\end{center}
\caption{Permeability and solution for the islands problem in 2d.}
\label{fig:islands_2d}
\end{figure}

The first test problem considers the scalar elliptic PDE \eqref{eq:ScalarBVP} in the unit square
or unit cube with the two-dimensional coefficent field given in Figure \ref{fig:islands_2d}. In the
three-dimensional version, the coefficient does not depend on the $z$-coordinate.
Boundary conditions are of Dirichlet type on the two planes perpendicular to the $x$-direction
and homogeneous Neumann on the rest of the boundary.

\subsubsection{Basic two-level experiments}

First, we gather some basic experiments that illustrate the behavior of spectral domain
decomposition methods. In \cite{Effendiev2010}, it was demonstrated that isolated large 
diffusion coefficients lead to very small eigenvalues that are well
separated from the rest. The spectrum
of the local GEVP for an interior subdomain contains zero eigenvalues
corresponding to the kernel
of the bilinear form~$a$ (i.e., the constant function here, or the
rigid body modes for linear elasticity), then a set of very small eigenvalues
related to isolated large coefficients and finally, with some
  gap, more or less equidistantly spaced eigenvalues.
In the spectral DD preconditioner, one has the choice of either
including a fixed number $n_{ev}$ of eigenvectors
per subdomain into the coarse space or to select all eigenvectors where the corresponding eigenvalue
is below a threshold $\eta$. In the first case, the size of the coarse space is controlled,
while in the second case the convergence rate is controlled.
In most experiments reported below, we will choose the basis for the
coarse space according to a threshold $\eta$.

\begin{table}[tbp]
\caption{Iteration numbers \#IT and coarse space sizes $n_{0}$ for a
  fixed number of subdomains in the two-level method, 
when varying fine mesh size $h$ and overlap $\delta$ (using
  Conjugate Gradients, 
$\mathbb{Q}_1$ conforming FEs in 2d, 16 subdomains and a fixed threshold $\eta=0.15$).}
\label{fig:basic1}
\begin{center}
\begin{tabular}{r|rrrr|rrrr}
\toprule
 & \multicolumn{4}{c|}{Laplace} &  \multicolumn{4}{c}{Islands}\\
 & \multicolumn{2}{c}{$\delta\sim h$} & \multicolumn{2}{c|}{$\delta\sim H$} 
 & \multicolumn{2}{c}{$\delta\sim h$} & \multicolumn{2}{c}{$\delta\sim H$} \\
$h^{-1}$ & \#IT & $n_{0}$ & \#IT & $n_{0}$ & \#IT & $n_{0}$ & \#IT & $n_{0}$ \\
\midrule
320   & 30 & 56 & 30 & 56    & 31 & 65   & 31 & 65\\ 
640   & 29 & 114 & 32 & 50  & 27 & 103 & 30 & 59 \\
1280 & 27 & 236 & 25 & 54  & 27 & 240 & 25  & 78\\
2560 & 26 & 488 & 26 & 54  & 25 & 481 & 24  & 55\\
\bottomrule
\end{tabular}
\end{center}
\end{table}
 Table \ref{fig:basic1} investigates how the size of the coarse space in the two-level method 
depends on the overlap $\delta$. For the standard two-level Schwarz method
it is well known that the convergence rate depends on $H/\delta$ \cite{TosWid05}. For fixed $H$
(fixed number of subdomains) and decreasing mesh size $h$ the number of iterations does not change
when $\delta\sim H$, while it will increase for $\delta\sim h$. For a fixed threshold $\eta$, Table 
\ref{fig:basic1} shows that here, the number of iterations $\#IT$ remains constant independent of the choice
of the overlap. However, the size of the coarse space $n_{0}$ increases when $\delta\sim h$ while it
does not increase when $\delta\sim H$. The table also shows that the size of the coarse space
increases only slightly when the homogeneous diffusion coefficient (Laplace) is changed into a heterogeneous 
diffusion coefficient (Islands).

\begin{table}[tbp]
\caption{Iteration numbers \#IT and coarse space sizes $n_{0}$ for a
  fixed number of subdomains in the two-level method, 
when varying polynomial degree $p$ and overlap $\delta$ (using Conjugate Gradients,  DG FEs
in 2d, $384^2$ elements and 256 subdomains).}
\label{fig:basic2}
\begin{center}
\begin{tabular}{rr|rr|rrr|rr|rrr}
\toprule
 & & \multicolumn{5}{c|}{$\eta=0.15$} & \multicolumn{5}{c}{$n_{ev}=20$} \\
 \midrule
 & & \multicolumn{2}{c|}{$\delta=2$} & \multicolumn{3}{c|}{$\delta$ var} & \multicolumn{2}{c|}{$\delta=2$} & \multicolumn{3}{c}{$\delta$ var} \\
$p$ & $n_1$ & \#IT & $n_{0}$ & $\delta$ & \#IT & $n_{0}$ & 
\#IT & $n_{0}$ & $\delta$ & \#IT & $n_{0}$\\
\midrule
1 & 589824  & 28 & 1457 & 2 & 28 & 1457& 18 & 5120& 2 &18 & 5120\\
2 & 1327104 & 21 &3171 & 3 &22 & 1901 & 19 & 5120& 3 & 18 &  5120\\
3 &2359296  & 20 & 5026 & 3 &21 & 2991 & 20 & 5120& 3 & 19 &  5120\\
4 & 3686400 & 18 & 8217 & 4 &21 & 3322 & 21 &5120 & 4 & 20 &  5120 \\
5 & 5308416 & 17 & 13596 & 4 &21  & 5078 & 23 &5120 & 4 & 21 &  5120 \\
6 & 7225344 & 17 & 17029 & 5 &22 & 5234 & 24 & 5120 & 5 & 22 & 5120  \\
\bottomrule
\end{tabular}
\end{center}
\end{table}

Table \ref{fig:basic2} investigates the two-level method applied to the DG discretization
of the Islands problem. The number of sudomains, as well as the mesh size is fixed in this computation
and the polynomial degree $p$ is varied. The number of degrees of freedom $n_1$ on the fine level is
increasing correspondingly. Experiments with a fixed threshold $\eta=0.15$ or a fixed number of eigenvectors
per subdomain $n_{ev}=20$, as well as with fixed and varying overlap $\delta$ are conducted.
Using a fixed threshold, the number of iterations is constant or even decreasing 
while the size of the coarse space $n_0$ increases with increasing polynomial
degree. With a fixed number of eigenvalues the iteration numbers are slightly increasing at a constant size of the coarse space.
In all cases, the preconditioner shows very good performance.

\subsubsection{Weak scaling in 2d}

\begin{table}[tbp]
\caption{Conjugate Gradients, $\mathbb{Q}_1$ conforming finite elements, $\delta=3h$, $\eta=0.3$.}
\label{fig:islands2d_weak_scaling}
\begin{center}
\begin{tabular}{l|rrrrr}
\toprule
subdomains & 64 & 256 & 1024 & 4096 & 16384 \\
\midrule
levels & \multicolumn{5}{c}{degrees of freedom}\\
\midrule
finest total & 410881 & 1640961 & 6558721 & 26224641 & 104878081\\
\midrule
2 lvl $n_0$ & 306 & 1348 & 5523 &  22673 & 91055\\
3 lvl $n_0$ &        & 130 & 431 & 1319 & 3890\\
4 lvl $n_0$ &        &        &  207 &  436  & 891 \\
\midrule
\midrule
levels & \multicolumn{5}{c}{iterations \#IT}\\
\midrule
2 & 25 & 26 & 27 & 26 & 26 \\
3 &      & 32 & 31 & 31 & 33\\
4 &      &      &  40 & 38 & 38\\ 
\bottomrule
\end{tabular}
\end{center}
\end{table}

We now turn to the multilevel method and carry out experiments with a varying number of subdomains.
Table \ref{fig:islands2d_weak_scaling} conducts a weak scaling experiment for the two-dimensional Islands problem,
where the number of degrees of freedom per subdomain is fixed. $\mathbb{Q}_1$ conforming finite elements
with a fixed overlap $\delta=3h$ and threshold $\eta=0.3$ are used. From left to right the number of subdomains
increases from 64 to 16384. The row labelled ``finest total'' gives the total number of degrees of freedom on the
finest level, while the next three rows report $n_0$, the size of the level 0 space when 2, 3 or four levels are used.
These results show that the size of the coarsest space can be significantly reduced (sizes of intermediate levels are not shown).
Finally, the last two rows give iteration numbers when using two,
three and four levels. Within each row we observe
robustness w.r.t.~the number of subdomains. Within each column
we observe a moderate increase with the number of levels,
but certainly not the exponential increase predicted by Theorem
\ref{thm:multilevelrate}. The numbers in Table
\ref{fig:islands2d_weak_scaling} would suggest $\kappa(BA) = O(L^2)$.

\subsubsection{Strong scaling in 3d}

Table \ref{fig:islands3d} gives results for the Islands problem in three dimensions using a cell-centered finite
volume discretization with two-point flux approximation. Here, the mesh is fixed and the number of  subdomains
as well as the number of levels are varied. The first set of rows corresponds to the two-level method. We note that the
sequential run-time $T_{seq}$ for setting up the preconditioner is reduced by almost a factor 3 when the 
number of subdomains is increased from 512 to 4096. 
This is due to the fact that the direct solver and the
eigensolver scale nonlinearly with the number of degrees of freedom
per subdomain, i.e., smaller is better. But at the same time the size of the coarse problem $n_0$ is increasing. The estimated parallel computation
time therefore has a minimum at 2048 subdomains with 25.2 seconds. With 4096 subdomains the
time for (sequential) factorization of the coarse problem $T_{coarse}$
becomes very large. Also note that there is quite a lot of
variability in the times needed to solve the eigenproblems in each subdomain. Minimum and maximum
times over all subdomains are reported in the columns labeled
$T_{i,min}$ and $T_{i,max}$, respectively. This suggests that
more subdomains than available processors should be used in
order to average runtimes over several subdomains.
The last three rows show corresponding results for a three level method using 4096 subdomains and different numbers
of subdomains on the intermediate level. The minimal parallel runtime
is achieved for 128 subdomains on level 1, leading to an
improvement over the two-level method in that case. Also note that in Table \ref{fig:islands3d} we only report times for constructing the preconditioner.
The solution time is
only $1/20$ of the setup time.
\begin{table}[tbp]
\caption{Islands problem in 3d. Fixed problem size $320^3$ mesh, 32768000 degrees of freedom.
On the finest level $n_{ev}=15$ eigenvectors are taken per subdomain, while for the three level calculation
threshold $\eta=0.4$ is used on the intermediate level. Iteration numbers are for the hybrid form of the preconditioner using
multiplicative subspace correction over levels and restricted additive Schwarz in each level used within GMRES (restart not reached).
Times are in seconds.}
\label{fig:islands3d}
\begin{center}
\begin{tabular}{rr|rrrrrrr}
\toprule
 $P_L$ & $P_{L-1}$ & \#IT & $n_0$ & $T_{seq}$ & $T_{par}$ & $T_{i,min}$ & $T_{i,max}$ & $T_{coarse}$\\
\midrule
\multicolumn{9}{c}{two level method}\\
\midrule
  512 & 1 & 12 & 7680  & 63613 & 191.3 & 70.3 & 176.2 & 0.47\\
1024 & 1 & 12 & 15360 & 35817 & 58.4 & 18.2 & 49.8 & 1.3\\
2048 & 1 & 14 & 30720 & 18781 & 25.2 & 4.9 & 13.2 & 5.1\\
4096 & 1 & 13 & 61441 & 19982 & 33.5 & 2.2 & 7.0 & 20.1\\
\midrule
\multicolumn{9}{c}{three level method}\\
\midrule
4096 &  32& 15 & 1387 & 21168 & 55.9 & 9.8 & 42.3 & 0.27 \\
4096 &  64& 15 & 1817 & 20725 & 27.7 & 2.5 & 15.1 &  0.18\\
4096 &  128&16 & 2569& 20549 & 18.4 & 0.59 &  6.2 &0.15 \\
\bottomrule
\end{tabular}
\end{center}
\end{table}

\subsection{SPE10 Problem}

\begin{figure}[tbp]
\begin{center}
\includegraphics[width=0.7\textwidth]{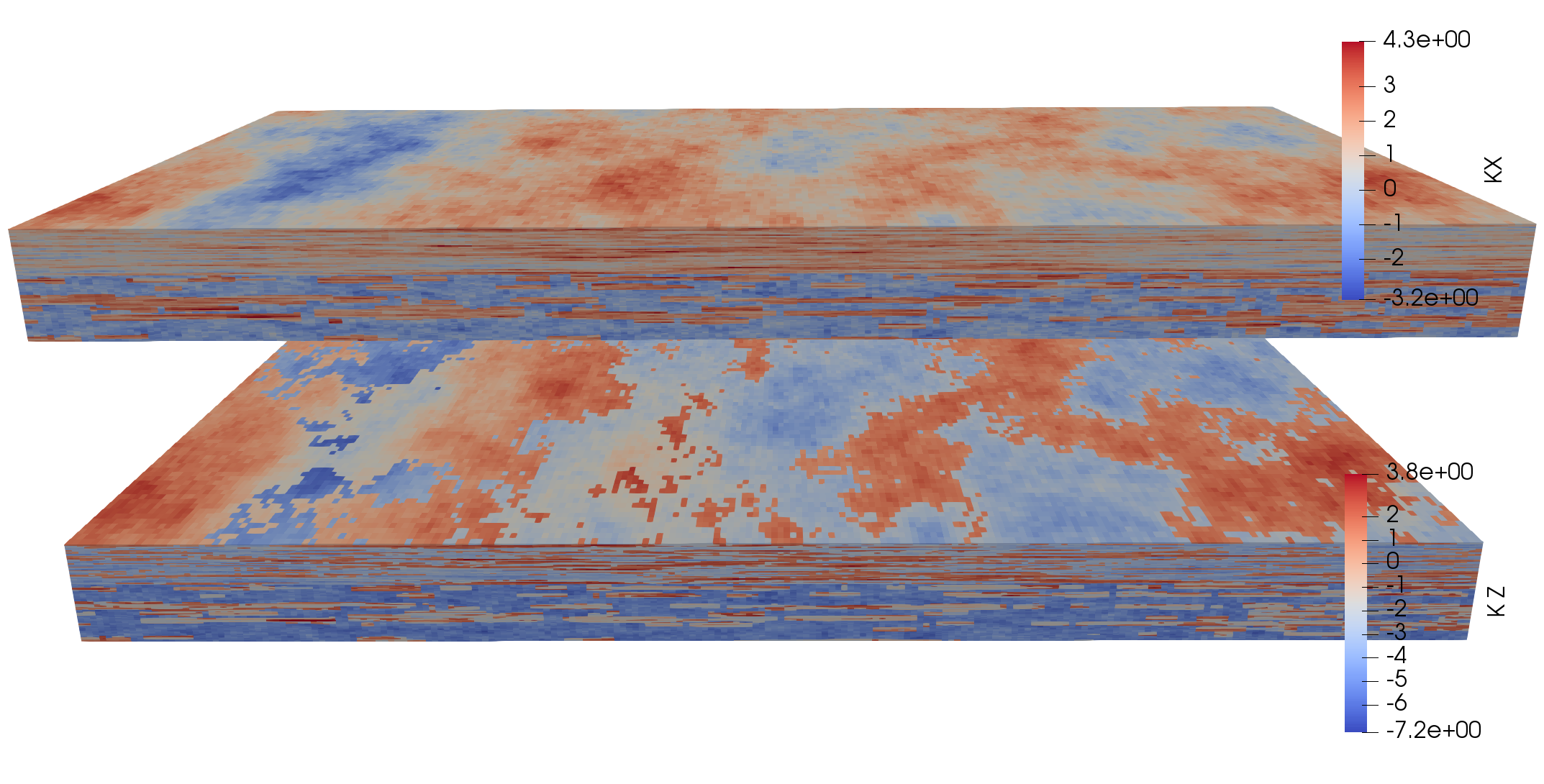}
\end{center}
\caption{Permeability field for the SPE10 problem.
}
\label{fig:spe10}
\end{figure}

\begin{table}[tbp]
\caption{SPE10 problem. Discretization schemes: 
$\mathbb{Q}_1$ conforming finite elements (CG), cell-centered finite volumes (FV) and
DG-$\mathbb{Q}_1$ are compared for two and three levels. 
Times are in seconds.
}
\label{fig:spe10_strong_scaling}
\begin{center}
\begin{tabular}{r|rrr|rrr|rrr}
\toprule
 & \multicolumn{3}{c|}{CG, $n_L=9124731$} & \multicolumn{3}{c|}{CCFV, $n_L=8976000$} & \multicolumn{3}{c}{DG, $n_L=8976000$} \\
$P_L, P_{L-1}$ & \#IT & $n_0$ & $T_{par}$ & \#IT & $n_0$ & $T_{par}$ & \#IT & $n_0$ & $T_{par}$ \\
\midrule
\multicolumn{10}{c}{two levels, $\eta = 0.3$}\\
\midrule
256 &   24  &  7237   &   133.8 &  25 & 7502   &  49.2   & 22    &  8366   &  233.5  \\
512 &  24   &  9830   &  60.4    & 23  & 10600 &   21.6  &  21   &  13690   & 124.9  \\
1024 & 28  &  21881   & 22.3   & 25&  25753&  12.3     &  24   &  31637   & 42.4  \\
2048 & 25  &  29023   &  15.7  &  25&  35411&  11.2    & 25    &  46844   & 36.7 \\
\midrule
\multicolumn{10}{c}{three levels, $\eta = 0.3$}\\
\midrule
256, 16&   29   & 1222     & 151.4      & 29     & 1364    & 70.3      &  31   &  1683   &  273.2  \\
512, 16&   27   &  1228    &  77.8    &  28    &  1446   &  47.4     &  28   &  1762   &  186.8  \\
1024, 32& 36   & 3145     &  46.3     &  34    & 3487    &  47.3     &  33   &  5476   &  231.4  \\
2048, 32& 31   & 3120     &  40.1     & 35     &  3421   &  49.9     &  36   &  5359   &  204.8  \\
\bottomrule
\end{tabular}
\end{center}
\end{table}

Next, we consider the SPE10 problem \cite{Christie2001}. Originally
intended as a benchmark for multiscale methods it is often used as a test problem for preconditioners as well.
The scalar elliptic problem \eqref{eq:ScalarBVP} is solved in a box-shaped domain, discretized with
an axiparallel and equidistant hexahedral mesh consisting of 1122000 elements. The diffusion tensor $K(x)$
is diagonal and highly variably. The $x$ and $y$ components are identical and vary over 7 orders of magnitude.
The $z$ component varies over 11 orders of magnitude. Figure \ref{fig:spe10} shows the permeability.

Table \ref{fig:spe10_strong_scaling} reports results for the SPE10 problem where we concentrate on the
comparison of the performance for different discretization schemes:
conforming $\mathbb{Q}_1$ finite elements
on a refined mesh, cell-centered finite volumes on a refined mesh and DG-$\mathbb{Q}_1$ on the original mesh.
All problems have roughly the same number of degrees of freedom, i.e., around 9 million. Results for two and three levels
using up to 2048 subdomains are given. The hybrid form of the preconditioner using
multiplicative subspace correction over levels and restricted additive Schwarz in each level is used within GMRES (restart not reached)
For each configuration we report number of iterations, size of the coarsest space and estimated parallel runtime for
setting up the preconditioner. We observe: the number of iterations is independent of the number of subdomains
and the discretization scheme used. From two to three levels a moderate increase in the number of iterations is observed.
However, the problem size is too small to achieve an improvement
in runtime through the use of more than two levels.

\subsection{Composites Problem}

We report results on modelling
carbon fibre composite materials from aerospace engineering, described in detail in
  \cite{REINARZ2018269,BUTLER2020106997}.
The setup is similar to the one in \cite[p. 271]{REINARZ2018269}, except
that the domain is flattened out and consists of only 9 ply layers and
8 interface (resin) layers. The equations of linear elasticity are solved 
in three dimensions using $\mathbb{Q}_2$ serendipity elements resulting in 10523067 degrees of freedom.
Table \ref{fig:composites} shows results for 1024 subdomains using 2,
3 or 4 levels using the preconditioner 
in its hybrid form within GMRES (multiplicative over levels,
restricted additive Schwarz within levels, restart not reached).
While the two-level method converges in 13 steps, the three and
  four level methods need 31 and 35 iterations, respectively. 
The maximum number of degrees of freedom in any coarse subdomain is
significantly reduced in the three and four level
methods compared to the two level method.

\begin{table}[tbp]
\caption{Carbon fibre composite problem. 9 ply layers (thickness 0.23mm),
8 resin layers (thickness 0.02mm), discretized with $256\times 64 \times 52$ mesh using
$\mathbb{Q}_2$ serendipity elements. Threshold $\eta=0.35$.}
\label{fig:composites}
\begin{center}
\begin{tabular}{rrrrr|rrrr|r}
\toprule
& \multicolumn{4}{c|}{subdomains} & \multicolumn{4}{c|}{max dofs/subdomain} & \#IT\\ 
\midrule
level &  3 & 2 & 1 & 0 & 3 & 2 & 1 & 0 \\
\midrule
&            & & 1024 & 1 &  & & 28791 & 21565 & 13\\
&        &1024 & 32 & 1 &  &28791 & 1260 & 914 & 31\\
&1024 & 128 & 16 & 1 & 28791 & 546& 515 & 273 & 35\\
\bottomrule
\end{tabular}
\end{center}
\end{table}

\section{Conclusions}
\label{sec:conclusions}

In this paper we extended the GenEO coarse space introduced in
\cite{Geneo14} from two to multiple levels and  
used it in the construction of multilevel preconditioners which are
robust in the fine mesh size, number of subdomains  
and coefficient variations.
The number of levels in the hierarchy is typically moderate since
aggressive coarsening is used. While the theory predicts an
exponential increase of the condition number of the preconditioned
system with the
number of levels numerical results suggest that the increase is
moderate. We believe that novel approximation theory for related
spectral coarse spaces in \cite{Chupeng} will allow
us to improve these theoretical results in future work.
In addition, the theory presented is more general than \cite{Geneo14},
extending also to different discretization schemes as well as to
different variants of the
generalized eigenproblem. In particular, we were able to analyse the
preconditioner for discontinuous Galerkin discretizations of scalar
elliptic problems. Numerical results illustrate
the robustness of  
the preconditioner for heterogeneous diffusion as well as linear
elasticity problems. Improvements over the two-level 
method could be demonstrated for a three-dimensional problem with 
30 million degrees of freedom. 

%
\appendix

\section*{Acknowledgments}

This work is supported by the Deutsche Forschungsgemeinschaft (DFG, German Research Foundation) 
under Germany's Excellence Strategy EXC 2181/1 - 390900948 (the Heidelberg STRUCTURES Excellence Cluster).
P. B. would like to thank Hussam Al Daas for discussions.

\bibliographystyle{siamplain}
\bibliography{lit}
\end{document}